\documentclass[a4paper,12pt,oneside]{amsart}

\usepackage{amsmath,amsfonts,amssymb,amsthm,amsopn,amsxtra}

\newtheorem{Thm}{Theorem}    
\newtheorem*{Cor}{Corollary}
\newtheorem{Corn}{Corollary}
\newtheorem{Lem}{Lemma}
\newtheorem{Pro}{Proposition}

\theoremstyle{remark}
\newtheorem{Rem}{Remark}

\newtheorem*{Di}{Discussion}
\newtheorem*{Ex}{Example}

\DeclareMathOperator{\supp}{supp}

\newcommand{\Cal}{\mathcal}
\newcommand{\fr}{\mathfrak}
\newcommand{\C}{\mathbb{C}}        
\newcommand{\R}{\mathbb{R}}        
\newcommand{\Rp}{\mathbb{R}^+}        
\newcommand{\SL}{\mathit{SL}}
\newcommand{\Or}{\mathit{O}}
\newcommand{\SO}{\mathit{SO}}
\newcommand{\SU}{\mathit{SU}}
\newcommand{\U}{\mathit{U}}

\newcommand{\starp}{\star_\R}
\newcommand{\odp}{\odot_\R}
\newcommand{\squ}{\mathbin{\square}}
\newcommand{\squp}{\mathbin{\square_\R\!}}

\begin{document}
\title[algebras of measures]{Arens products\\[2mm]
for some convolution algebras of measures.} 
\author{Viktor  Losert}
\address{Institut f\"ur Mathematik, Universit\"at Wien, Strudlhofg.\ 4,
  A 1090 Wien, Austria}
\email{Viktor.Losert@UNIVIE.AC.AT}
\date{2 December 2025}

\begin{abstract}
We consider the measure algebra of a Ch\'ebli-\!Trim\`eche hypergroup
(in particular, double coset spaces of classical Lie groups) and study
the corresponding Arens products on its second dual. The behaviour turns
out to be different to the group case investigated in \cite{LN}. For this,
we study more closely properties of the multiplication and generalized
translation in a Ch\'ebli-\!Trim\`eche hypergroup and the asymptotic
behaviour.
\end{abstract}
\maketitle
\vspace*{-.7cm}

\baselineskip=1.3\normalbaselineskip
\setcounter{section}{-1}
\section{Introduction}\vspace{1mm} \label{Intro}
Let $M$ be a Banach algebra with multiplication $\star$ and dual space $M'$.
Taking up the notations of  \cite{LN}, we define actions on $M'$ and $M''$
\begin{eqnarray*}
\langle h\cdot\mu,\,\nu\,\rangle & = & \langle\,h\,,\mu\star\nu\rangle
	\text{\;\;for\;\;} h\in M',\;\mu,\nu\in M\,,\\
\langle\fr{n}\odot h,\,\mu\,\rangle & =
& \langle\,\fr{n}\,,h\cdot\mu\rangle \text{\;\;\;for\;\;}
    \fr{n}\in M'',\; h\in M',\; \mu\in M\,,\\
\langle\fr{m}\squ\fr{n},\,h\,\rangle & = &
    \langle\,\fr{m}\, , \fr{n} \odot h \rangle \text{\;\;for\;\;}
    \fr{m},\fr{n}\in M'',\; h\in M' .
\end{eqnarray*}
$\squ$ is called the \emph{left} (or first) \emph{Arens product} on
$M''$, it makes $M''$ a Banach algebra. $\odot$ makes $M'$ a left
$M''$-module. The
(left) {\it topological centre} of $M''$ is defined as 
\[ Z_t(M'')\ = 
\{\; \fr{m} \in M''  : \;
\text{the\;mapping\ \;} \fr n\mapsto\fr m\squ\fr n
\text{\ is\;weak* continuous\;on\;} M'' \,\}\,.
\]
A second method to extend the
multiplication of $M$ \,(interchanging the order of the factors above)
gives the right (or second) Arens product. Always
$M\subseteq Z_t(M'')$ and
$M$ is called {\it strongly Arens irregular} if $Z_t(M'')=M$ \;and the
same holds for the right topological centre. See
\cite{Da}\;p.\,248ff. for a further discussion of Arens products and
topological centres and \cite{DaLa}. If $M$ is abelian, $Z_t(M'')$
coincides with the centre of $(M'',\squ)$ in the algebraic sense and also
with the right topological centre of $M''$.

Let $G$ be a locally compact group. In \cite{LL} it was shown that
$L^1(G)$ (with convolution $\star$) is strongly Arens irregular.
Then in \cite{LN}, the same was proved for the measure algebra
$M(G)$ \,(space of bounded real or complex Radon measures).\vspace{1mm}

Let $K$ be a compact subgroup of $G$\,,
then $G/K=\{xK:x\in G\}$ denotes the space of left cosets,
$G//K=\{KxK:x\in G\}$ the space double cosets (with quotient topology).
$G/K$ and $G//K$ are locally compact spaces and it is well known (see also
\cite{LN}\,Lemma\,1\,) that
there is a canonical isomorphism of the Banach space $M(G/K)$ with
$M(G)\star\lambda_K$ \;($\lambda_K$\,: normalized Haar measure of $K$)
and similarly there is an isomorphism of $M(G//K)$
with $\lambda_K\star M(G)\star\lambda_K$\,. Thus, modulo these isomorphisms,
we may consider as subalgebras \,$M(G//K)\subseteq M(G/K)\subseteq M(G)$, in
particular, there is defined a generalized convolution for measures on
$G//K$ resp. $G/K$\,.

When $K$ is normal in $G$\,, $G//K=G/K$ is a group.
For a general compact subgroup $K$\,, by \cite{BH}\,p.\,55, $H=G//K$ is a
{\it hypergroup}.
For a general hypergroup~$H$ there is given a mapping $\star$ assigning to
$x,y\in H$ a probability measure
$\delta_x\star\delta_y$  on $H$ (we write $\delta_x$
for the point measure at $x$) and some axioms
are satisfied, see \cite{BH}\,1.1.2. This extends to a {\it generalized
convolution} for
bounded measures on $H$\,, defining
$\mu\star\nu=\iint\delta_x\star\delta_y\,d\mu(x)\,d\nu(y)$ for
$\mu,\nu\in M(H)$
\,(\cite{BH}\,p.\,10, they use the notation $M^b(\cdot)$ for our
$M(\cdot)$\,). In this way $M(H)$ becomes a Banach algebra (for any hypergroup
$H$). For a bounded Borel measurable function $f\!:H\to\C$\,, {\it generalized
translations} are defined by
$(T_yf)(x)=\langle\delta_x\star\delta_y\,,\,f\,\rangle$
for $x,y\in H$ \,(\cite{BH}\,1.2.1). See also \cite{Z} for some aspects of
hypergroups.
\begin{Ex}\ \ $G=\SL(2,\C),\ K=\SU(2)$
\\
Put $\alpha(t)=\begin{pmatrix}e^{\frac t2} & 0\\0& e^{-\frac t2}
\end{pmatrix}$. With $\Rp=[0,\infty[$\,, it follows from elementary linear
algebra (singular values) that $\{\alpha(t)\!:t\in\Rp\}$ constitutes a
complete set of representatives for the family of double cosets $G//K$\,.
This gives rise to a homeomorphism of $\Rp$ and $G//K$\, which induces an
isometric (w*-continuous) isomorphism of the spaces of measures $M(\Rp)$ and
$M(G//K)$.
Combined with the embedding of $M(G//K)$ into $M(G)$ mentioned above,
$\delta_t\in M(\Rp)$ is mapped to
$\mu_t=\lambda_K\star\delta_{\alpha(t)}\star\lambda_K\in M(G)\ \,(t\ge0)$.
\\[.2mm]
For \;$x,y>0$ \;consider \
$k_{xy}(t)=\begin{cases}
\dfrac{\sinh(t)}{2\,\sinh(x)\,\sinh(y)}\ &\text{for \ }
\lvert x-y\rvert<t<x+y \\[2.8mm]
\hspace{1.cm} 0&\text{else}\hspace{1.5cm} (t\ge0). 
\end{cases}$
\\[1mm]
Computing $\mu_x\star\mu_y$ and transforming back to $\Rp$ gives
\;$\delta_x\star\delta_y=k_{xy}\,\lambda$\quad
where $\lambda$ \,denotes standard Lebesgue
measure on $\Rp$ \,(see \cite{Jew}\,(15.5D), \,in his notation
$d_t=\alpha(2t)$, another factor~2 comes silently through his mapping
$\varphi$ \,-- see 15.5B; \;a corresponding formula is given in
\cite{BH}\,p.237). Furthermore, $\delta_x\star\delta_0=\delta_0\star\delta_x=
\delta_x$\,. $H=\Rp$\,, equipped with this convolution is called the
{\it Naimark hypergroup} (\cite{BH} 3.5.66).
\end{Ex}
The Naimark hypergroup belongs to the class of
{\it Ch\'ebli-\!Trim\`eche hypergroups} \,(see the Appendix for some remarks
on their definition). In Theorem\,1 we determine  $Z_t(L^1(H)'')$  when $H$
is a Ch\'ebli-\!Trim\`eche hypergroup. It has an additional (``remote") part
$Z^\infty(H)$. Corollaries\;1,2 describe cases where $Z^\infty(H)$ is trivial
(hence $L^1(H)$ strongly Arens irregular) or non-trivial
(hence $L^1(H)$ not strongly Arens irregular). Theorem\,2 is devoted to the case
of $M(H)$. There $Z_t(M(H)'')$ has another component $M_0(H)''_{\;LA}$ which
is always non-trivial (Proposition\,3). Thus
$M(H)$ is {\it not} strongly Arens
irregular for all Ch\'ebli-\!Trim\`eche hypergroups $H$\,. For $G//K$\,, this
covers the cases where $(G,K)$
is a Gelfand pair of non-compact type and real rank $1$. More explicitly
(\cite{BH}\,p.\,246), this includes $\SL(2,\C)$ as above, also
$G=\SL(2,\R),\,K=\SO(2)$ \,(and more generally, $G$ a connected semisimple Lie
group with finite centre and real rank $1$, $K$ a maximal compact subgroup).
Here $L^1(G//K)$ is strongly Arens 
irregular (Corollary\,1).
Furthermore, there are the Euclidian motion groups 
$G=\R^n\rtimes K,\,K=\SO(n)$ (or $K=\Or(n)$\,),
$n\ge2$ and here $L^1(G//K)$ turned out to be {\it not} strongly Arens
irregular (Corollary\,2). Note that $\R^2\rtimes\SO(2)$ is solvable,
hence this type of behaviour is not restricted to non-amenable groups.
Related cases for compact groups and hypergroups are described in Remark\,4.
In Remark\,6 we correct some wrong assertions in the literature.
\\[2mm]
$Z^\infty(H)$ and $M_0(H)''_{\;LA}$ are defined by annihilating properties for
products. For the proof of the results we investigate more deeply (Section\,2)
properties of the product and translations in
Ch\'ebli-\!Trim\`eche hypergroups (Proposition\,1) that are quite different
from the group case, using results of Braaksma, de Snoo (\cite{BS}) and
Stein, Wainger (\cite{SW}). Asymptotically there are
relations to standard convolution on $\R$\,, depending on properties of the
Sturm--Liouville function $A$ (Proposition\,2), obtained using results of
Brandolini, Gigante (\cite{BG})\,.
\\[3mm]
{\bf Some basics.} \ For $f\in M',\;\mu\in M$,
the value of the functional is written as $\langle f,\mu\rangle$.
As usual, $M$ is identified with a subspace of the second dual
$M''$ by the canonical embedding (which amounts to
$\langle \mu,f\rangle=\langle f,\mu\rangle$). If $N$ is a subspace of $M$\,,
then $N''$ is identified with the w*-closure of $N$ in $M''$ (equivalently, the
image of the second dual of the inclusion mapping). If $N$ is a subalgebra
of $M$\,, it is easy to see that in this way $N''$ becomes a subalgebra of
$M''$ and if $N$ is an ideal in $M$\,, then $N''$ is an ideal in $M''$
and the left action $\odot$ of $N''$ on $N'$ extends to an action (also
denoted by $\odot$) of $M''$ \;(in particular,
$\fr m\odot(\fr n\odot h)=(\fr m\squ\fr n)\odot h$
for $\fr m,\fr n\in M'',\,h\in N'$).
For a subset $W\subseteq M\,,\ W^\perp\subseteq M'$ will denote its
annihilator (the functionals vanishing on $M$).\vfill
\section{Main results}\vspace{2mm} \label{Main}
In \cite{BH} (and further references given there) quite a number of
results can be found showing analogies to locally compact abelian groups,
description of the dual, Plancherel theorem, etc.
We collect now some properties of the convolution in a
Ch\'ebli-\!Trim\`eche hypergroup, exhibiting significant differences to the
group case (and providing the technical basis for the results on the
Arens products), see also Section\,\ref{Proofs},\,Remark\,8.
Recall that the underlying space of $H$ is always
$\Rp=[0,\infty[$\,, differences come from the convolution $\star$\,, which is
defined by a Sturm--Liouville function $A$ \,(see the Appendix). To simplify
the proofs we will assume that $A$ is $C^2$ on $]0,\infty[$\,.
$\lambda$~denotes {\it standard} Lebesgue measure on $\Rp$ (or $\R$)
and we put $L^1(H)=L^1(\Rp,\lambda)$. We write
$\starp$ for standard convolution on the group $\R$\,.\vspace{-3mm}
\begin{Pro}	
Let $H$ be a Ch\'ebli-\!Trim\`eche hypergroup.
\item[(a)]\ For $x,y>0,\ \,\delta_x\star\delta_y$ is absolutely continuous
with respect to $\lambda$ (shortly:\linebreak $\delta_x\star\delta_y\in L^1(H)$).
$(x,y)\mapsto\delta_x\star\delta_y$ is \,{\rm norm}--continuous
(on $]0,\infty[\times]0,\infty[$)\,.\vspace{1mm}
\item[(b)]\ For $x>0,\quad
\nu_x=\lim\limits_{y\to\infty}\,\delta_{-y}\!\starp(\delta_x\star\delta_y)$
exists in the norm topology of $L^1(\R)$, convergence is uniform for $x$ on
compact subsets of $]0,\infty[$\,. \;$x\mapsto\nu_x$ is norm--continuous.
\;For $x\to\infty$, we have $\nu_x\to0$ in the weak* topology (i.e.
$\sigma(M(\R),\linebreak C_0(\R))$\,).
\item[(c)]\ Let $f\!:\Rp\to\C$ be bounded, Lebesgue-measurable.
Then $T_xf$ is continuous on $]0,\infty[$ for $x>0$.
$\lim\limits_{x\to0+}T_xf=f$ \,holds pointwise a.e. on $\Rp$.
\end{Pro}\vspace{-2mm}
Additional results on the asymptotic behaviour of $(\nu_x)$ for $x\to\infty$
can be given under conditions on the Sturm Liouville function $A$\,. We
restrict to major cases covering the examples of double coset spaces mentioned
before.\vspace{-1mm}
\begin{Pro}    
Let $H$ be a Ch\'ebli-\!Trim\`eche hypergroup.
\item[(a)]\ If $H$ has sub-exponential growth and $A$ is unbounded then
\\
$\lim_{y\to\infty}\lVert\,\delta_x\starp\nu_y-\nu_y\rVert=
\lim_{y\to\infty}\lVert\,\nu_{y+x}-\nu_y\rVert=0$ \;holds for
every $x\in\Rp$ \;(asymptotic invariance).
\item[(b)]\ If \,$0<\lim_{y\to\infty}A(y)\,e^{-2\rho y}<\infty$\,, with
$\rho>0$,
then \,$\nu_\infty=\lim_{y\to\infty}\delta_{-y}\!\starp\nu_y$ \,exists in
the norm topology and the Fourier transform of \,$\nu_\infty\in L^1(\R)$
satisfies $\hat\nu_\infty(t)\neq0$ for all $t\in\R$\,.\vspace{-2mm}
\end{Pro}\vspace{-2mm}
\noindent The proofs of Propositions\,1 and 2 will be given in 
Section\;\ref{Proofs}. See also Remarks\,\,8,\,9 for further\vspace{2mm} comments.

We state some consequences for the Banach algebras $M(H)$. \pagebreak Recall that
$M(H)$ is {\it commutative} (\cite{BH}\,Cor.\,3.4.3) and $\delta_0$ is the
identity element. Put
\;$M_0(H)=\linebreak
\{\mu\in M(\Rp)\!:\, \mu(\{0\})=0\}\quad  (\,\cong\;M(]0,\infty[)\,)$.
\begin{Cor}[\bf to Proposition\,1]\ {\rm (a)}	    
\ For \,$0<a<b<\infty$ \ the set\newline
$\{\,\mu\star\nu:\;\mu,\nu\in M([a,b]), \
\lVert\mu\rVert,\!\lVert\nu\rVert\le1\}$ \;is norm--compact in $L^1(H)$.
\vspace{1mm}
\item[(b)]\ $M_0(H)\star M_0(H)\subseteq L^1(H)$.\vspace{-1mm}
\end{Cor}
\noindent Thus $M_0(H)$ is a subalgebra of $M(H)$, \;$M(H)$ is (isomorphic to)
the unitization of~$M_0(H)$\,, \;$L^1(H)$ is an ideal in $M(H)$, $M_0(H)/L^1(H)$
has trivial multiplication.
\begin{proof}
By Proposition\,1(a), $\{\delta_x\star\delta_y:a\le x,y\le b\}$ is
norm-compact in $L^1(H)$, hence by Mazur's theorem, the same is true for its
norm-closed absolutely-convex hull. Then this is w*-compact (in $M(H)$) as
well, in particular w*-closed. But, by the definition of convolution
(\cite{BH}\,p.\,10), for $\mu,\nu\in M([a,b])$ with
$\lVert\mu\rVert,\!\lVert\nu\rVert\le1$,\linebreak
$\mu\star\nu$ belongs to
the w*-closed absolutely-convex hull of
$\{\delta_x\star\delta_y:a\le x,y\le b\}$ \;(in particular
$\mu\star\nu\in L^1(H)$). It follows from \cite{BH}\,Th.\,1.2.31 (or
\cite{Jew}\,L.\,2.4B) that $(\mu,\nu)\mapsto\mu\star\nu$ is w*-continuous on
bounded, uniformly tight subsets of $M(H)\times M(H)$, hence in our case,
it is \,w*\,to\;norm continuous on the unit ball of $M([a,b])\times M([a,b])$.
This proves (a). Now (b) follows easily by approximation. \ (Note that in
\cite{BH}\,Th.\,1.2.31 the proof is given only for non-negative measures. The
more general assertion there is not true in general and needs additional
assumptions, similarly in \cite{BH}\,Prop.\,1.2.16(iv)\,).
\end{proof}
\begin{Rem}	
\ It is an easy consequence to describe the Gelfand
duals (character spaces): \;$M_0(H)\sphat=L^1(H)\sphat$
\ (the non-zero multiplicative
functionals on $L^1(H)$ have unique extensions to $M_0(H)$\,),
$M(H)\sphat=L^1(H)\sphat\ \cup\{\varphi_0\}$ with
$\langle\varphi_0,\mu\rangle=
\mu(\{0\})$ \,(see also \cite{CMS}\,Th.\,4.6 for a related case).
$L^1(H)\sphat$ \;can be described in terms of the functions $\phi_\lambda$
appearing in section\;\ref{Proofs}, proof of Proposition\,1(b) \;(\cite{BH}\,Th.\,3.5.50).
\end{Rem}
\begin{Di} {\rm ($\alpha$) }
In some sense one can ``localize" certain elements of $L^1(H)''$ as follows.
If $B$ is a Borel subset of $\Rp$ (with $\lambda(B)>0)$ there is a canonical embedding of
$L^1(B)$ into $L^1(H)$ (extending functions by $0$) and a canonical
projection (restricting functions). The second duals of these mappings
give an embedding of $L^1(B)''$ into $L^1(H)''$ and a projection (both w*-continuous). In this way,
$L^1(B)$ will be considered as a subspace of $L^1(H)$, identifying it with
$L^1(H)_B=\{u\in L^1(H): u=0 \text{ outside } B\}$. $L^1(B)''$ will be
identified with the w*-closure of $L^1(H)_B$ in $L^1(H)''$.
We say that $\fr m\in L^1(H)''$ is {\it carried} by $B$ if it belongs to
$L^1(B)''$
(compare \cite{LP}\,p.\,449) and that $\fr m$ {\it vanishes} on $B$ if
its projection to $L^1(B)''$ is $0$ \,(i.e. $\fr m$ is carried by $\Rp\setminus B$\,). The
same for $M(H)''$. Alternatively, the algebra
$B(\Rp)$ of bounded Borel measurable functions with pointwise multiplication
acts on the space $M(H)$, making also $M(H)''$ a $B(\Rp)$--module 
(acting by w*-continuous operators on $M(H)''$) and then the projection
is given by \,$c_B\,\fr m$ where
$c_B$ denotes the indicator function of $B$\,.
\item[($\beta$)] Analogously to ($\alpha$)\;, $L^\infty(B)$ will be considered as a subspace of $L^\infty(H)$.
As in \cite{LP}\,p.\,448, we write $L_0^\infty(H)$ for the norm closure of
\ $\bigcup_{b<\infty}\!L^\infty([0,b])$  \;(equiv\-alently, these are the
elements $f\!\in\!L^\infty(H)$ such that
$\lVert fc_{]b,\infty[}\rVert_\infty\!\to\!0$ for
$b\!\to\!\infty$).\linebreak Its dual $L_0^\infty(H)'$ can be identified with
the norm closure of \ $\bigcup_{b<\infty}\!L^1([0,b])''$ and\linebreak
its annihilator is given by
\,$L_0^\infty(H)^\perp=\{\,\fr m\!\in\!L^1(H)''\!:
\langle \fr m,f\rangle=0 \text{ \ for all}\linebreak
f\!\in\! L^\infty([0,b]),\;b<\infty\,\}$, i.e. $\fr m$ vanishes on all
bounded subsets of $\Rp$ ($\fr m$ is  ``located at infinity"). As in
\cite{LP}\,Th.\,2.8,
$L_0^\infty(H)^\perp$ is a w*-closed ideal, $L_0^\infty(H)'$ a norm closed
subalgebra and \,$L^1(H)''=L_0^\infty(H)'\oplus L_0^\infty(H)^\perp$.
Furthermore, $L_0^\infty(H)^\perp$ is also a left ideal in
$(L^1(\R)'',\squp)$\,.
\item[($\gamma$)] Put $S\mu=\int\nu_x\,d\mu(x)$ \,(with $\nu_x$ from Proposition\,1(b)). This
defines a linear contractive map \,$S\!:M_0(H)\to L^1(\R)$
\,(recall that $x\mapsto\nu_x$ is continuous on $]0,\infty[$\,).
Write \,$S_1=S\vert L^1(H)$ \;($S_1\!\!:L^1(H)\to L^1(\R)$\,).
For $f\in L^\infty(\R,\lambda),\ x>0$ put
\;$(Tf)(x)=\langle\nu_x,f\rangle$\,. Then
\,$T\!:L^\infty(\R,\lambda)\to C_b(]0,\infty[)\subseteq L^\infty(\Rp,\lambda)$
\,is a contractive linear operator. It is easy
to see that $T=S_1'$ \,(dual mapping), consequently $T$ is w*-continuous.
\item[($\delta$)] The left action of $M(\R)''$ on $L^1(\R)'$ arising from $\starp$ is written
as \,$\odp$\,. Using the composition~$\starp\ ,\ L^1(\Rp)$ becomes a
subalgebra of $L^1(\R)$ \,(but $L^1(\Rp)'$ is not
$L^1(\Rp)$-invariant under
the action $\odp$).\vspace{.5mm}
\end{Di}

The asymptotic behaviour of $\delta_x\star\delta_y$ when $x$ or $y$ is
large will now be used to give some relationship between $\squ$ and $\squp$ and also between
$\odot$ and $\odp$\,, involving the operators $S,T$ introduced in ($\gamma$).\vspace{-1mm}

\begin{Lem}\ {\rm (a) }    
For $\fr n\in L_0^\infty(H)^\perp\,,\
\fr m\in M_0(H)''\;(\supseteq\!L^1(H)''\,)$
we have \,$\fr m\squ\fr n=S''\fr m\,\squp\,\fr n$\,,
and for $f\in L^\infty(H)$ we have \;$\fr n\odot f=T(\fr n\odp f)$\,.
\\
If $\fr m\in M(H)''$ is carried by a compact subset of \,$]0,\infty[$\,,
then for $\fr n\in L_0^\infty(H)^\perp$ we have
\,$\fr n\squ\fr m=\fr n\,\squp S''\fr m=
\fr m\squ\fr n$\,.
\item[(b)]\ If \,$\fr m\in L^1(H)''$ is carried by $[0,b]$ for every $b>0$,
then $\fr n\squ\fr m=\langle\fr m,\mathbf{1}\rangle\,\fr n$ holds for every
$\fr n\in L^1(H)''$
and \,$\fr m\odot f=\langle\fr m,{\bf 1}\rangle\,f$ for $f\in L^\infty(H)$.
If \,$f\in B(\Rp)$ (bounded Borel measurable function) is continuous at $0$,
then
\,$\langle\fr m,f\rangle=\langle\fr m,{\bf 1}\rangle\,f(0)$\,.\vspace{-2mm}
\end{Lem}\noindent
$\mathbf{1}$ denotes the constant function with value $1$\,. 
Observe that in this Lemma (and the further arguments prior to Lemma\,3)
$\odot$ always refers to the action of $M(H)''$ on $L^\infty(H)$ (see the end
of the Introduction). For the treatment of $Z_t(M(H)'')$ an extended version of
Lemma\,1 (proved in the same way) will be given immediately before Lemma\,3, involving
some additional notation.\vspace{-2mm}
\begin{proof}
For (a), we consider $\mu\in M([a,b])$ for some
$0\!<\!a\!<\!b<\!\infty$\,. Given $\varepsilon>0$\,,\linebreak it follows from
Proposition\,1(b) that there exists $c>0$ such that
$\lVert\delta_x\star\delta_y-\linebreak\delta_y\starp\nu_x\rVert_1<\varepsilon$
for all $x\in[a,b]\,,\; y>c$\,. Observe that 
\,$\iint v(y)\,\delta_y\starp\nu_x\,d\mu(x)\,dy=\linebreak
\int v(y)\,S\mu\starp\delta_y\,dy= S\mu\starp v$\,. Hence
$\lVert\,\mu\star v-S\mu\starp v\,\rVert<
\varepsilon\,\lVert \mu\rVert\lVert v\rVert$ \,holds for $v\in L^1([c,\infty[)$.
\\
For $\fr n\in L_0^\infty(H)^\perp$ this gives by approximation
$\lVert \mu\squ\fr n-S\mu\squp\fr n\rVert\le
\varepsilon\,\lVert\mu\rVert\lVert\fr n\rVert$\,.
With $\varepsilon\to0$, we get $\mu\squ\fr n=S\mu\squp\fr n$ for
every $\mu\in M_0(H)$ and the first formula in (a) follows by approximation
of $\fr m$\,.
\\
For $\fr m\in M([a,b])''$ the formula above gives
$\lVert\fr m\squ v-S''\fr m\squp v\rVert\le
\varepsilon\,\lVert\fr m\rVert\lVert v\rVert$
\,for $v\in L^1([c,\infty[)$. The expression on the left equals
$ v\squ\fr m-v\squp S''\fr m$ and as above, we get
$\fr n\squ\fr m=\fr n\squp S''\fr m$ for
$\fr n\in L_0^\infty(H)^\perp$. Observe that by Proposition\,1(a) (as in the
proof of the Corollary\,(a) to Proposition\,1\,) the operator $S\vert M([a,b])$
is compact. Hence $S''\fr m\in L^1([a,b])$ and this implies that
$\fr n\squp S''\fr m=S''\fr m\squp\fr n$, consequently
$\fr n\squ\fr m=\fr m\squ\fr n$\,.
\\
The second formula in (a) follows easily from the first one by duality.
\\[.5mm]
For (b), we consider $\fr n=u\in L^1(H)$. For $\varepsilon>0$ there exists
$b>0$ such that
$\lVert u\star\delta_x-u\rVert<\varepsilon$ for $x\in[0,b]$
(see \cite{BH}\,1.4.6(ii), ``translation continuous measures").
Then for $v\in L^1([0,c])$ we get
$\lVert u\star v-\langle v,{\bf 1}\rangle\,u\rVert<
\varepsilon\,\lVert u\rVert\lVert v\rVert$\,. This implies the first formula
in (b) and the second one follows again by duality. If $f$ is continuous at
$0$ put $f_0=f(0)\,{\bf 1}$\,. For $\varepsilon>0$ there exists
$b>0$ such that $\sup_{t\in[0,b]}\lvert f(t)-f_0(t)\rvert<\varepsilon$. Since
$\fr m$ is carried by $[0,b]$ it follows that
$\lvert\langle\,\fr m\,,f-f_0\rangle\rvert<\varepsilon\,\lVert\fr m\rVert$\,.
This implies the last formula in (b).\vspace{-1.9mm}
\end{proof}
\begin{Rem}	
\ Be aware that our space $L^1(H)$ (based on standard Lebesgue
measure~$\lambda$) is different from the usage in
\cite{BH} and \cite{Jew} who work with the Haar measure $\lambda_A=A\lambda$
of $H$ \,(I believe, the present version makes the asymptotics at infinity
easier).
The algebra structure $\star$ for our $L^1(H)$ is defined using the embedding
$f\mapsto f\lambda$ into $M(H)$\,. Clearly $f\mapsto\frac fA$ defines an
isometric linear isomorphism  $L^1(\Rp,\lambda)\to L^1(\Rp,\lambda_A)$. 
Writing $\star_{BH}$ for the convolution of functions as defined in
\cite{BH}\,1.4.8 (corresponding to $f\mapsto f\lambda_A$), it is related
to our notation by $f\!\star_{BH}g=\frac1A((fA)\star(gA))$ \,(similarly for
$\mu\star f$ when $\mu\in M(H)$\,). Thus both $L^1$-spaces give realizations
of the same subalgebra $M_a(H)$ of $M(H)$ \,(\cite{BH}\,1.6.21).
\end{Rem}
\begin{Rem}	
\ It follows from Lemma\,1(b) that if
$\fr m\in\bigcap_{b>0}L^1([0,b])''$ (i.e. $\fr m$ is carried by all
neighbourhoods of $0$) and $\langle\fr m,\mathbf{1}\rangle=1$ then $\fr m$
is a right identity in $L^1(H)''$. The argument carries over also to the
group case. Similar to \cite{LP}\,p.\,447, let $\Cal E$ be the set of right
identities in $L^1(H)''$, \,$\Cal E_1=\{\fr e\in\Cal E:\lVert\fr e\rVert=1\}$.
Then one can show that $\Cal E_1\subseteq\bigcap_{b>0}L^1([0,b])''$ and
$\bigcap_{b>0}L^1([0,b])''$ is the linear subspace generated by $\Cal E_1$\,.
Also $\Cal E=\Cal E_1+L^1(H)''_{\;RA}$\;(right annihilator) and
$L^1(H)''_{\;RA}=(L^\infty(H)\cdot L^1(H))^\perp$. Observe that in the
notation of \cite{BH} $f\cdot u=\int u(x)\,T^x\!f\,dx$, hence
$L^\infty(H)\cdot L^1(H)$ consists of the bounded $\beta$-uniformly
continuous functions on $H$ in the sense of \cite{BH}\,Def.\,1.2.26(ii). 
\end{Rem}\vspace{-.5mm}
As in \cite{Da}\,Def.\,3.3.29 we write
$L^1_0(H)=\{u\in L^1(H):\langle u,\mathbf{1}\rangle=0\}$
for the {\it augmentation ideal} (functions with integral $0$).\vspace{-2.4mm}
\begin{Lem}\ {\rm (a) }    
$L^1_0(H)$ is the norm closure of the linear space generated by\linebreak
$\{\delta_x\star u-u:\ x\in H\,,\,u\in L^1(H)\}$\,.
\item[(b)]\ $L^1_0(H)$ has a bounded approximate identity.\vspace{-1.5mm}
\end{Lem}\noindent
This is the analogue of results for a locally compact group $G$\,.
(a) corresponds to \cite{RS}\,Ch.3,\;6.4 (using the notation $J^1(G,G)$).
By \cite{Da}\,Th.\,5.6.42, existence of a bounded approximate identity for
$L^1_0(G)$ is equivalent to amenability of $G$\,.
\begin{proof}
For (a) we use duality. The annihilator of the set consists of
$f\in L^\infty(H)$ such that
$\langle f,u\star\delta_x\rangle=\langle f,u\rangle$ for all
$x\in H\,,\,u\in L^1(H)$\,. By approximation, it follows that
$\langle f,u\star v\rangle=\langle f,u\rangle\,\langle v,\mathbf{1}\rangle$
for all $v\in L^1(H)$\,. Taking an approximate identity $(u_n)$ it follows
that 
$\langle f,v\rangle=\lim\,\langle f,u_n\rangle\,\langle v,\mathbf{1}\rangle$
\,thus $f=\lim\langle f,u_n\rangle\,\mathbf{1}$ must be constant.
\\
For (b) take $v_n=\frac1n\,c_{[0,n]}$\,. For $x\in H$\,, estimates as in the
proof of Proposition\,1(a) above give
\,$\lim_{n\to\infty}\,\lVert \delta_x\star v_n-\nu_x\starp v_n\rVert=0$ \;and
by elementary arguments\linebreak
$\lim\,\lVert\nu_x\starp v_n-v_n\rVert=0$\,. Thus
\;$\lim\,\lVert\delta_x\star v_n-v_n\rVert=0$ for all $x\in H$ and then it
follows
from (a) that $\lim\,\lVert u\star v_n\rVert=0$ for all $u\in L^1_0(H)$.
If $(u_n)$ is an approximate identity for $L^1(H)$ as in the proof of
Theorem\,1 below, then $(u_n-v_n)$ is a bounded approximate identity for
$L^1_0(H)$.
\end{proof}\vspace*{-2mm}
\begin{Thm}	    
Let $H$ be a Ch\'ebli-\!Trim\`eche hypergroup.
Put
\\
$Z^\infty(H)=\{\fr m\in L^1(H)\squ L_0^\infty(H)^\perp:\
\fr m\squ L_0^\infty(H)^\perp=\{0\}\,,
\ \fr m\odot L^\infty(H)\subseteq C_0(H)\,\}$.
\\[.5mm]
Then \ $Z_t(\,L^1(H)'')=L^1(H)\oplus Z^\infty(H)$\,.
\end{Thm}
\begin{proof}
In one direction, $L^1(H)\subseteq Z_t(\,L^1(H)'')$ is true in general. Take
$\fr m\in Z^\infty(H)$. Then
$\fr m\odot L^\infty(H)\subseteq C_0(H)\subseteq L_0^\infty(H)$ implies
$\fr n\squ\fr m=0=\fr m\squ\fr n$ for $\fr n\in L_0^\infty(H)^\perp$.
\\[.5mm]
If $\fr n\in L^1(H)''$ is carried by a compact subset of $]0,\infty[$\,\,,
then by Lemma\,1(a) \,$\fr m\squ\fr n=\fr n\squ\fr m$\,.
\\
Let $(u_n)_{n\ge1}\subseteq L^1(H)$ be such that
$u_n\ge0\,,\,\lVert u_n\rVert=1$ and
$u_n$ vanishes outside $[0,\frac 1n]$\,. Then $(u_n)$ is an approximate
identity for $L^1(H)$ (\cite{BH}\,Rem.\,1.6.16). From\linebreak
$\fr m\in L^1(H)\squ L_0^\infty(H)^\perp$ it follows that
\,$\lim u_n\squ\fr m=\fr m$\,.
Then for $f\in L^\infty(H)\linebreak \fr m\odot f\in C_0(H)$ implies \
$\fr m\odot f\,(0)=\lim\,\langle\fr m\odot f,u_n\rangle=
\lim\,\langle u_n\squ\fr m,f\rangle=\langle\fr m,f\rangle$\,.
If $\fr n\in L^1(H)''$ is carried by $[0,b]$ for every $b>0$, this gives
using Lemma\,1(b)
$\langle\fr n\squ\fr m,f\rangle=\langle\,\fr n\,,\fr m\odot f\rangle=
\langle\fr n,\mathbf{1}\rangle\,\fr m\odot f\,(0)=
\langle n,\mathbf{1}\rangle\langle\fr m,f\rangle=
\langle\,\fr m\,,\fr n\odot f\rangle=\langle\fr m\squ\fr n,f\rangle$\,.
\\[.4mm plus .5mm]
It is easy to see (see also ($\beta$) below) that every
$\fr n\in L_0^\infty(H)'$ splits into $\fr n=\fr n_0+\fr n_1$
where $\fr n_0$ is carried by $[0,b]$ for every $b>0$ and $\fr n_1$ belongs to
the norm closure of $\bigcup\{L^1([a,b]):0<a<b<\infty\}$\,. It follows
from the special cases above that again $\fr n\squ\fr m=0=\fr m\squ\fr n$\,,
hence $\fr m\in Z_t(\,L^1(H)'')$.
\\
$Z^\infty(H)\subseteq L_0^\infty(H)^\perp$ implies
$L^1(H)\cap Z^\infty(H)=\{0\}$\,, giving a direct sum.
\\[.9mm plus .5mm]
For the converse take now $\fr m\in Z_t(\,L^1(H)'')$. In three steps we will
treat parts of $\fr m$ arising from projections to various locations.
\item[($\alpha$)]
Let $(u_n)$ be an approximate identity as above. It is well known that
every w*--accumulation point $\fr e$ of $(u_n)$ in $L^1(H)''$ is
a right identity (i.e. $\fr n\squ\fr e=\fr n$ for all $\fr n\in L^1(H)''$,
see e.g.\,\cite{LP}\,p.\,447; this is also a special case of Lemma\,1(b)\,).
Centrality of $\fr m$ implies \,w*--$\lim\,\fr m\squ u_n=\fr m$\,.
\\
Let $\fr m_{a,b}$ be our projection of $\fr m$
to $L^1([a,b])''$ (see Discussion\;($\alpha$)\,). We will show now that
\,$\fr m_{a,b}\!\in L^1(H)$ for all \,$0<a<b<\infty$\,. Recall that
\,\,$\supp(\delta_x\star\delta_y)\subseteq
[\mspace{1.5mu}\lvert x-y\rvert\mspace{1.5mu},x+y\;]$ holds
for all $x,y\in\Rp$ (\cite{BH}\,Th.\,3.5.36). Take some $\varepsilon$ with
$0<\varepsilon<a$ \,and put $a'=a-\varepsilon\,,\,b'=b+\varepsilon$\,. If $u\in L^1(H)$
is carried by $[0,\varepsilon]$ and $v$ is carried by
$[0,a']$, then $v\star u$ vanishes on $[a,b]$,
the same if $v$ is carried by $[b',\infty[$\,. It follows
that $(\fr m\squ u)_{a,b}=(\fr m_{a',b'}\squ u)_{a,b}$\,.
\\
For an $u'\in L^1(H)$ that is carried by a compact subset of \,$]0,\infty[$
\,we can apply the Corollary to Proposition\,1, arguing similarly as in \cite{Ld}\;Prop.\,1.
Recall that on a norm--compact subset the w*-topology coincides with the norm topology.
Take a bounded net $(v_\iota)$ in $L^1(H)$, w*-converging to $\fr m_{a',b'}$  with $v_\iota$
carried by $[a',b']$. It follows that $(u'\star v_\iota)$ is norm-convergent, hence
\,\,$u'\squ \fr m_{a',b'}=\lim\,u'\star v_\iota\in L^1(H)$. The set of these $u'$ is norm--dense
in $L^1(H)$, by approximation we arrive at \,$\fr m_{a',b'}\squ u\in L^1(H)$ for all
$u\in L^1(H)$. Then also $(\fr m_{a',b'}\squ u)_{a,b}\in L^1(H)$\,.
\\
Now we can use weak sequential completeness.
\\
Just for the record: \;the technique invoking weak sequential completeness
in this context
was already used in \cite{LL} (see L.\,5) to manage the case of metrizable
groups, then also in \cite{LL2}.
Then it was taken up and elaborated in more detail in \cite{LU},
without reference, and some people spoke of a "new method"\dots
\\
In the present case, we consider the approximate identity $(u_n)$ above. Then (the projection is
w*--continuous) 
w*--$\lim\,(\fr m\squ u_n)_{a,b}=\fr m_{a,b}$ and
$(\fr m\squ u_n)_{a,b}=\linebreak(\fr m_{a',b'}\squ u)_{a,b}\in L^1(H)$ holds for $n>\frac 1a$\;. 
The sequence $(\,(\fr m\squ u_n)_{a,b})_{n>\frac1a}$
\,is weakly Cauchy in $L^1(H)$ and weak sequential completeness of $L^1(H)$
implies $\fr m_{a,b}\in L^1(H)$\,.\vspace{1mm plus 2mm}
\item[($\beta$)] Combining the pieces obtained from ($\alpha$) \,(varying $a,b$), we get
$w\in L^1(H)$ such that $\fr m_{a,b}=w_{a,b}$ for all $0<a<b<\infty$\,.
Replacing $\fr m$ by $\fr m-w$ we may assume now that $\fr m_{a,b}=0$ for all
$0<a<b<\infty$\,.
Then $\fr m_0=\fr m_{0,b}$ does not depend on $b>0$, hence
$\fr m_0$ is carried by $[0,b]$ for every $b>0$\,. We will show next that
$\fr m_0=0$.
\\[2mm]
For $u\in L^1(H)$ Lemma\,1(b) implies
$\fr m_0\squ u=\langle\fr m_0,{\bf 1}\rangle\,u$\,. If $u$ is
carried by $[0,\delta]$, where $0<\delta<b$ it follows as in ($\alpha$)
that
\,$(\fr m\squ u)_{0,b}=(\fr m_{0,b+\delta}\squ u)_{0,b}=(\fr m_0\squ u)_{0,b}=
\langle\fr m_0,{\bf 1}\rangle\,u$\,. For the approximate identity this implies
w*--$\lim\,\langle\fr m_0,{\bf 1}\rangle\,u_n=\fr m_0$\,. If
$\langle\fr m_0,{\bf 1}\rangle$ would
be non--zero, then $(u_n)$ would be weakly Cauchy, hence as above, weakly
convergent. The limit would be given by a function vanishing for $x>0$ which
is impossible. Thus $\langle\fr m_0,{\bf 1}\rangle=0$ and using the limit
again it follows that $\fr m_0=0$. \;(Or more directly: we might choose $u_n$
carried by $[\frac1{n+1},\frac1n]$ and $f(x)=1$ on
$[\frac1{n+1},\frac1n]$ for $n$ even and $f(x)=0$ otherwise, then it is
immediate that $\langle u_n,f\rangle$ diverges).\vspace{0mm plus 1mm}
\item[($\gamma$)] We have now that $\fr m_{0,b}=0$ for all $b>0$, i.e.
$\fr m\in L_0^\infty(H)^\perp$. We want to show that $\fr m\in Z^\infty(H)$
which will finish the proof.
\\
By Lemma\,1(a) and centrality, the mapping $\fr n\mapsto S_1''\fr m\squp\fr n$
is w*-continuous for $\fr n\in L_0^\infty(H)^\perp$. We use now Neufang's
factorization technique: \;there exists $f_0\in L^\infty(\R)$ such that
$\{\,\fr n\odp f_0:\;\fr n\in L_0^\infty(H)^\perp,\ \lVert\fr n\rVert\le1\}$
contains the unit ball of $L^\infty(\R)$ \,(see \cite{N}\,Th.\,2.1). Then, as
in the proof of \cite{LN}\,Lemma\,3
it follows that $S_1''\fr m\in L^1(\R)$. Now take $f\in C_0(\R)$ and
$\varepsilon>0$. By Proposition\,1(b) there exists $c>0$ such that
$\lvert\langle \nu_t,f\rangle\rvert<\varepsilon$ for $t>c$\,. Since
$\fr m\in L^1([c,\infty[)''$, it follows (approxima\-ting by a net from
$L^1([c,\infty[)$\,) that
$\lvert\langle S_1''\fr m,f\rangle\rvert<\varepsilon\,\lVert\fr m\rVert$\,.
This gives\linebreak $\langle S_1''\fr m,f\rangle=0$\, and then
$S_1''\fr m\in L^1(\R)$ implies $S_1''\fr m=0$\,. Consequently, by
Lemma\,1(a), \,$\fr m\squ L_0^\infty(H)^\perp=\{0\}$\,.
\\
It follows that for $f\in L^\infty(H)$ we have
$\fr m\odot f\in L_0^\infty(H)^\perp$. Centrality of $\fr m$ implies that
$\langle\,\fr n\,,\fr m\odot f\rangle=\langle\fr m\squ\fr n,f\rangle$ \,for
\,$\fr n\!\in\!L^1(H)''$.
Hence if \,$\fr m=\text{w*-}\lim\fr m_\alpha$ \,it follows that $\fr m\odot f$
is the weak limit of $(\fr m_\alpha\odot f)$ \,\,(i.e.\;for
$\sigma(L^\infty(H),L^1(H)'')$\,).
For $\fr m_\alpha\!\in L^1(H)$ one has \,$\fr m_\alpha\odot f\in C_b(H)$
\,\,(e.g.\,by \cite{Jew}\,Th.\,5.5D), hence \,$\fr m\odot f\in C_b(H)$\,. This
implies that $\fr m\odot f\in C_0(H)$\,.
\\
For an approximate identity $(u_n)$ as above,
\,w*--$\lim\,\fr m\squ u_n=\fr m$
was noted in ($\alpha$). Now the Grothendieck property of
$L^\infty(H)$ implies that this is a weak limit. Hence $\fr m$ belongs to the
weak closure of \,\,$L^1(H)\squ L_0^\infty(H)^\perp$, but by the factorization
theorem, this is already closed.
\end{proof}
The size of $Z^\infty(H)$ depends on the asymptotic behaviour of the
measures $\nu_x$ from Proposition\,1(b).
\begin{Corn}   
If \,$\lim_{y\to\infty}\lVert\delta_x\starp\nu_y-\nu_y\rVert=
\lim_{y\to\infty}\lVert\nu_{y+x}-\nu_y\rVert=0$ \,holds for
every $x\in\Rp$, one has \,$Z^\infty(H)=L^1_0(H)\squ L_0^\infty(H)^\perp=
L^1_0(\R)\squp L_0^\infty(H)^\perp$.
\\
Hence $L^1(H)$ is not strongly Arens irregular.\vspace{-2mm}
\end{Corn}\noindent
By Proposition\,2(a) this covers the case of the Bessel--Kingman
hypergroups (\cite{BH} 3.5.61), in particular the
double coset hypergroups coming from the Euclidian motion groups, as
mentioned in the Introduction.\\
It is not hard to see that
$L^1_0(H)\squp L_0^\infty(H)^\perp$ is always non--trivial.
For example, let $(x_n)\subseteq H$ be such that $x_{n+1}-x_n\,\to\infty$
and \,$\fr m$ a w*-accumulation point of $(\delta_{x_n})\subseteq M(H)''$.
Then $\mu\mapsto\mu\squ\fr m$ is isometric on $M(H)$ and
\,$L^1(H)\squ\fr m\subseteq L_0^\infty(H)^\perp$\,.
\begin{proof}
In fact, $Z^\infty(H)\subseteq L^1_0(H)\squ L_0^\infty(H)^\perp$ holds
for any Ch\'ebli-\!Trim\`eche hypergroup $H$\,: \ Take
$v_n\in L^1([n,\infty[)$ with $v_n\ge0\,,\,\lVert v_n\rVert=1$\,. For
$\fr m\in Z^\infty(H)\,,\linebreak \fr m\odot L^\infty(H)\subseteq C_0(H)$
implies
w*--$\lim v_n\odot\fr m=0$\,. Again by the Grothendieck property it follows
that this tends weakly to $0$\,. If $(u_n)$ is chosen as before, we get
that $(u_n-v_n)\squ\fr m\to\fr m$ weakly, thus $\fr m$ belongs to the weak
closure of $L^1_0(H)\squ L_0^\infty(H)^\perp$\,.
\\[1mm]
Since $T\mathbf{1}=\mathbf{1}$\,, the subspace $L^1_0(H)$ is $S_1$--invariant. Hence
Lemma\,1(a) implies $L^1_0(H)\squ L_0^\infty(H)^\perp\subseteq
L^1_0(\R)\squp L_0^\infty(H)^\perp$\,.
\\
For the converse, take $\fr m\in L^1_0(\R)\squp L_0^\infty(H)^\perp$ and we
will show that $\fr m\in Z^\infty(H)$\,. Using the analogue of Lemma\,2(a)
for $G=\R$ (\cite{RS}\,Ch.3,\,\,6.4), it is enough to consider
$\fr m=\delta_x\squp\fr m_1-\fr m_1$ with
$\fr m_1\in L^1(\R)\squp L_0^\infty(H)^\perp,\ x\in \Rp$.
If $(u_n)$ is chosen as before, observe that (since $\nu_x$ is
carried by $[-x,x]$\,) \;$(S_1u_n)$ is an approximate identity for $L^1(\R)$.
It follows from Lemma\,1(a) that $\fr m\in L^1(H)\squ L_0^\infty(H)^\perp$\,.
\linebreak
$\lim_{y\to\infty}\lVert\nu_{y+x}-\nu_y\rVert=0$ implies (making similar
estimates as before) that\linebreak $S_1''(\delta_x\squp\fr m_1-\fr m_1)=0$\,,
thus $\fr m\squ L_0^\infty(H)^\perp=\{0\}$.
$\lim_{y\to\infty}\lVert\delta_x\starp\nu_y-\nu_y\rVert=0$ gives for
$f_1\in L^\infty(\R)$ that $T(\delta_x\odp f_1-f_1)\in L_0^\infty(H)$,
hence by Lemma\,1(a) we have for $f\in L^\infty(H)\,,\
\fr m\odot f\in L_0^\infty(H)$\,. Furthermore, $\fr m=\lim u_n\squ\fr m$
implies that
$\fr m\odot f=\lim u_n\odot(\fr m\odot f)$ is continuous at $0$\,. Combined
this gives $\fr m\odot L^\infty(H)\subseteq C_0(H)$\,.
\end{proof}
\begin{Corn}	  
If \,$\nu_\infty=\lim_{y\to\infty}\delta_{-y}\starp\nu_y\in L^1(\R)$ \,exists
in the norm topology and its Fourier transform
satisfies $\hat\nu_\infty(t)\neq0$ for all $t\in\R$, one has
\,$Z^\infty(H)=\{0\}$\,.
\\
Hence $L^1(H)$ is strongly Arens irregular.\vspace{-2mm}
\end{Corn}\noindent
By Proposition\,2(b) this covers the case of the Jacobi
hypergroups of non-compact type (\cite{BH}\,3.5.64), in particular the
double coset hypergroups coming from the semisimple groups desribed
in the Introduction, e.g. the Naimark hypergroup.
\begin{proof}
Making similar estimates as before, the relation
$\nu_\infty=\lim_{y\to\infty}\delta_{-y}\starp\nu_y$ implies that
$S_1''\fr m=\fr m\squp\nu_\infty$ for $\fr m\in L_0^\infty(H)^\perp$. For
$\fr m\in Z^\infty(H)$ the property $\fr m\squ L_0^\infty(H)^\perp=\{0\}$ is
equivalent to $S_1''\fr m=0$ \,(use Lemma\,1(a) and recall that
$L_0^\infty(H)^\perp\odp L^\infty(\R)=L^\infty(\R)$\,, e.g. by
\cite{N}\,Th.\,2.1). Hence in our case $\fr m\squp\nu_\infty=0$. By
Wiener's theorem (\cite{RS}\,Ch.1,\,4.1\,and\;4.3), $\nu_\infty$ generates
a dense ideal in $L^1(\R)$, thus $\fr m\squp u=0$ for all $u\in L^1(\R)$.
Taking an approximate identity in $L^1(H)$, we get
$\fr m=\lim\,\fr m\squp S_1u_n=0$\,.
\end{proof}
\begin{Corn}	
If $\fr m\in Z_t(\,L^1(H)'')$ and $\fr m\ge0$ then $\fr m\in L^1(H)$.
\\
If $\fr m\in Z_t(\,L^1(H)'')$ is carried by $[0,b]$ (for some $b>0$) then
$\fr m\in L^1(H)$.\vspace{5mm}
\end{Corn}
Now we pass to $M(H)''$.
As an analogue of $L_0^\infty(H)^\perp$ and $L_0^\infty(H)'$, we have a
decomposition $M_0(H)''=M_{\rm fin}\oplus M_{\rm inf}$\,,
putting $M_{\rm inf}=\bigcap_{c>0}M([c,\infty[)''$ and $M_{\rm fin}$ shall
be the norm closure of \,\,$\bigcup_{c>0}M(]0,c])''$.
It is easy to see that $M_{\rm fin}\cap L^1(H)''=L_0^\infty(H)'$ and
$M_{\rm inf}\cap L^1(H)''=L_0^\infty(H)^\perp$. To get more
explicit examples, we consider $M_d(]0,\infty[)\ (\cong l^1(]0,\infty[)\,)$
the subspace of discrete (purely atomic) measures in
$M_0(H)$ \;(i.e., the norm-closed subspace generated by the point measures).
\\[1mm]
We have mentioned before that $\squ$ behaves well with respect to
subalgebras. But concerning the first duals and the left action $\odot$
more care is needed. To be more accurate, we write
$j_d\!:M_d(]0,\infty[)\to M_0(H)$ for the inclusion
mapping, similarly $j_a\!:L^1(H)\to M_0(H)$. Under duality,
$j_d'\!:M_0(H)'\to M_d(]0,\infty[)'\ (\cong l^\infty(]0,\infty[)\,)$ and
$j_a'\!:M_0(H)'\to L^1(H)'\ (\cong L^\infty(H)\,)$ are quotient mappings
(restricting the functionals). $B_\lambda(]0,\infty[)$ denotes the space of
bounded, {\it Lebesgue}
measurable functions \,$g\!:\;]0,\infty[\to\C$\,, considered as a subspace of
$M_d(]0,\infty[)'$, it contains $B(]0,\infty[)$
\,(bounded Borel measurable functions),
$C_b(]0,\infty[)$ and $C_b([0,\infty[)$ \,(bounded continuous functions).
For $g\in B_\lambda(]0,\infty[)\,,\; %
q(g)$ shall denote its equivalence class $\mod\lambda$ \,(thus
$q\!:B_\lambda(]0,\infty[)\to L^\infty(H)$\,).
For $g\in B(]0,\infty[)$, a functional $j_B(g)$ on $M_0(H)$ is defined
by $\langle j_B(g),\mu\rangle=\int\! g(x)\,d\mu(x)$ \,(often we are writing
briefly $\langle g,\mu\rangle$). Then
$j_B\!:B(]0,\infty[)\to M_0(H)'$ \,is an isometry, it extends the canonical
embedding $C_0(]0,\infty[)\to C_0(]0,\infty[)''=M_0(H)'$
(giving a more accessible
part of the dual $M_0(H)'$; similarly, if $g\in B([0,\infty[)$ one gets
$j_B(g)\in M(H)'$, in \cite{DS}\,p.\,44 this is denoted as
$\kappa_E(g)$\,). We have \,$j_d'(j_B(g))=g$ and \,$j_a'(j_B(g))=q(g)$\,.
\\
$\odot_0$ will now denote the action of $M(H)''$ on $M_0(H)'$ and
$\odot_a$ that of $M(H)''$ on $L^\infty(H)$ \;(formerly, in the investigation
of $Z_t(\,L^1(H)'')$, $\odot_a$ was just written as $\odot$).
For $\fr m\in M_0(H)'',\;\varphi\in M_0(H)'$ it follows from Corollary\,(b) to
Proposition\,1 that $\fr m\odot_0\varphi$ depends only on the
projection (restriction) $f=j_a(\varphi)$\,, hence we may write
$\fr m\odot_0 f$
for this. We have $j_a'(\fr m\odot_0 f)=\fr m\odot_a f$\,, thus for  
$\fr m\in M_0(H)'',\,f\in L^\infty(H),\ \;\fr m\odot_0 f$ is an extension of
the functional $\fr m\odot_a f$ on $L^1(H)$ to $M_0(H)$\,.
\\
In the same way as before, one can show now an {\it extended version}
of Lemma\,1(a)\,: \ 
For $\fr n\in M_{\rm inf}\,,\
\fr m\in M_0(H)''\;(\supseteq\!L^1(H)''\,)$\,,
we have \,$\fr m\squ\fr n=S''\fr m\,\squp\,\fr n$\,.
\\[1mm]
Also Lemma\,1(b) extends: \
Assume that \,$\fr m\in M_0(H)''$ is carried by $]0,b]$ for every $b>0$.
Then $\fr n\squ\fr m=\langle\fr m,\mathbf{1}\rangle\,\fr n$ holds for every
$\fr n\in M_0(H)''$
and \,$\fr m\odot f=\langle\fr m,j_B(\mathbf{1})\rangle\,f$ for
$f\in L^\infty(H)$.
If \,$g\in B(\Rp)$ is continuous at $0$,
then
\,$\langle\fr m,j_B(g)\rangle=\langle\fr m,j_B(\mathbf{1})\rangle\,g(0)$\,.
\\[1mm]
The formula $\fr n\odot f=T(\fr n\odp f)$ in Lemma\,1(a) was meant
(in the more precise notation introduced above) as
\;$\fr n\odot_af=q\circ T(\fr n\odp f)$ \,(where
$T\!:L^\infty(\R)\to C_b(]0,\infty[)$ \,and
$\fr n\in L_0^\infty(H)^\perp\subseteq L^1(H)''$). Also \;$S_1'=q\circ T$.
In the first
part of the proof of Lemma\,3(c) below, we will show that in fact
$j_d'(\fr m\odot_0f)=T(\fr m\odp f)$ and $\fr m\odot_0f=j_B(T(\fr m\odp f))$
for $\fr m\in M_{\rm inf}\,,\ f\in L^\infty(H)$\,.\vspace{-1.4mm}
\begin{Lem}\ {\rm (a) }    
For $\mu\in M_0(H)\,,\ f\!\in L^\infty(H)$ we have
\,$j_d'(\mu\odot_0f)\in C_b(]0,\infty[)$\,,
and \,$\mu\odot_0f=j_B(j_d'(\mu\odot_0f))$\,. If $g\in C_b([0,\infty[)$
then \,$j_d'(\mu\odot_0q(g))\in C_b([0,\infty[)$\,.
\item[(b)]\ If $\fr m\in M_{\rm inf}\,,\ f\in L^\infty(H)$\,, there exists
$g\in C_b([0,\infty[)$, such that
\,$j_d'(\fr m\odot_0f)=j_d'(\fr m\odot_0q(g))$\,.
\item[(c)]\ If \;$\fr m\in M_0(H)'',\ f\!\in L^\infty(H)$\,, then
\,$j_d'(\fr m\odot_0 f)\in B_\lambda(]0,\infty[)$ \,and
\,$q\circ j_d'\,(\fr m\odot_0 f)=\fr m\odot_a f$\,.
\end{Lem}
\begin{proof}
For $\nu\in M_0(H)$ we have (by approximation)
\,$\iint\langle\delta_x\star\delta_y,f\rangle\,d\nu(x)\,d\mu(y)=
\langle\nu\star\mu,f\rangle$\,. For
$\nu=\delta_x$ this describes $j_d'(\mu\odot_0f)(x)$ and continuity on
$]0,\infty[$ follows from Proposition\,1(a). The second formula in (a) follows
immediately from this. If $g$ is continuous on $H$ then
$(x,y)\to\langle\delta_x\star\delta_y,g\rangle$ is continuous on $H\times H$
by axiom HG3 for a hypergroup (\cite{BH}\,1.1.2).
\\[0mm plus.5mm]
For (b) note that
$j_d'(\fr m\odot_0f)(x)=\langle\,\fr m\,,\delta_x\odot_0 f\rangle$.
Since $\fr m\in M_{\rm inf}$ it will be enough to construct
$g\in C_b([0,\infty[)$
such that for every $x>0$ one has $
\langle\,\delta_x\star\delta_y\,,f-q(g)\rangle\to0$
for $y\to\infty$\,. This will be done by induction. Put $f_n=f\,c_{[n,n+1]}$
and choose any $g_0\in C([0,1])$. Assume that $g_{n-1}$ has been constructed.
Then, using compactness from the Corollary\,(a) to Proposition\,1, one can
approximate $f_n$ and adapt in a small neighbourhood of $n$ to find
$g_n\in C([n,n+1])$ \,(extended to $\Rp$ by $0$) such that
$\lVert g_n\rVert_\infty\le\lVert f\rVert_\infty$\,, $g_n(n)=g_{n-1}(n)$ and
\,$\lvert\langle\delta_x\star\delta_y\,,f_n-q(g_n)\rangle\rvert<\frac1n$
\,holds for \linebreak
$x\in[\frac1n,n]\,,\ y\in[n,n+1]$\,. Then putting $g=g_n$ on $[n,n+1]$
will do the job.
\\[1mm plus1mm]
For (c) assume first that $\fr m\in M_{\rm inf}$\,. For $\nu\in M_0(H)$ the
extended version of Lemma\,1(a) gives
\;$\nu\squ\fr m=S\nu\squp\fr m=\int\!\nu_x\!\squp\fr m\:d\nu(x)$\,.
This implies that
\linebreak
$j_d'(\fr m\odot_0f)=T(\fr m\odp f)$ which is continuous on $]0,\infty[$
\,(Proposition\,1(b)\,). In addition, the formula for $\nu\squ\fr m$ shows
that $\fr m\odot_0f=j_B(T(\fr m\odp f))$\,. Using the properties
of $j_a'\,,j_d'$\,, we get the last assertion of (c) in this case.
\\[0mm plus1mm]
Next we consider $\fr m\in M([a,b])''$ with $\lVert\fr m\rVert\le1$ where
$0<a<b<\infty$\,. Let $\mu\in M([a,b])$ be the restriction of $\fr m$ to
$C_0([a,b])\ (\subseteq M([a,b])')$. If $\fr m=\lim \mu_\alpha$ for $\sigma(M'',M')$ with
$(\mu_\alpha)\subseteq M([a,b])$, then $\mu=\lim \mu_\alpha$ for
$\sigma(M,C_0)$ (where $M=M_0(H),\;C_0=C_0(H)$). For any
$0<a_0\le a<b\le b_0<\infty$ we have the norm--compact set $C=
\{\,\rho\star\nu:\;\rho,\nu\in M([a_0,b_0]), \
\lVert\rho\rVert,\!\lVert\nu\rVert\le1\}\ (\subseteq L^1(H))$ from the
Corollary to
Proposition\,1. On $C$ both w*-topologies (\,$\sigma(M,C_0)$ and
$\sigma({\smash{L^1}}'',L^\infty)$\,)
coincide with the norm--topology. Hence the continuity properties of
$\star$ (see the proof of the Corollary to Proposition\,1) and $\squ$
imply that $\fr m\squ\nu=\mu\star\nu$ for all $\nu\in M([a_0,b_0])$. This
gives $\fr m\odot_0f=\mu\odot_0f$ and in that case (c) follows from (a).
\\[1mm plus1mm]
Finally assume that $\fr m\in\bigcap_{b>0}M(]0,b])''$. For $f\in L^\infty(H)$
observe that $j_d'(\delta_x\odot_0 f)=T_xf$\,. Put
$L_f=\{y>0:f_0(y)=\lim_{x\to0+}(T_xf)(y)\text{ \;exists}\,\}$ (some sort
of Lebesgue set for $f$). By Proposition\,1(c)
\;$\lambda(\,\Rp\!\setminus L_f)=0$ \,and (taking a sequence tending to $0$)
$f_0$ is Borel measurable on $L_f$\,. For $y\in L_f$ one gets from the
extension of Lemma\,1(b) that
$j_d'(\fr m\odot_0 f)(y)=\langle\fr m,j_B(\mathbf{1})\rangle\,f_0(y)$ which
implies that $j_d'(\fr m\odot_0 f)$ is Lebesgue measurable. In the same way,
for $v\in L^1(H)$ and a bounded sequence $(\mu_n)$ with
$\mu_n\in M(]0,\frac1n])\,,\ \langle \mu_n,\mathbf{1}\rangle=1$, dominated
convergence gives $\langle\mu_n\odot_0 f,v\rangle\to\langle f_0,v\rangle$\,,
from which one can deduce the second assertion of (c).
\\[0mm]
Combined one obtains (c) for all $\fr m\in M_0(H)''$.
\end{proof}\vspace{1mm}
The {\it left annihilator} of $M_0(H)''$ \;(now in the sense of algebras\,!\,)
is denoted as
\linebreak
$M_0(H)''_{\;LA}=\{\fr m\in M_0(H)'': \;\fr m\squ M_0(H)''=(0)\}$
\quad (\cite{Da}\,Def.\,1.4.4).\vspace{-1mm}
\begin{Pro}    
Let $H$ be a Ch\'ebli-\!Trim\`eche hypergroup. Then $M_0(H)''_{\;LA}$
is non--trivial.\vspace{-2.5mm}
\end{Pro}
\begin{proof}
Under duality, $\fr m\in M_0(H)''_{\;LA}$ is equivalent to
$\langle\,\fr m\,,\fr n\odot_0f\rangle=0$ for all
$\fr n\in M_0(H)'',\ f\in L^\infty(H)$\,.
If  $\fr m_d\in M_d(]0,\infty[)''$ vanishes on $B_\lambda(]0,\infty[)$, then
it follows from Lemma\,3(c) that \,$\fr m=j_d''(\fr m_d)\in M_0(H)''_{\;LA}$\,.
More generally, take any $\fr m_a\in L^1(H)''$ and let
$\fr m_d\in M_d(]0,\infty[)''$
be any extension of $\fr m_a\circ q$\,, then
\,$j_d''(\fr m_d)-j_a''(\fr m_a)\in M_0(H)''_{\;LA}$\,. Observe
(see also \cite{DS}\,p.\,28) that \
$j_d''\oplus j_a'':\linebreak
M_d(]0,\infty[)''\oplus L^1(H)''\;(l^1\text{-sum})\to M_0(H)''$
\ is isometric.
\end{proof}
\begin{Thm}	       
Let $H$ be a Ch\'ebli-\!Trim\`eche hypergroup.\vspace{-1mm}
\[Z_t(\,M(H)'')=M(H)\oplus Z^\infty(H)\oplus M_0(H)''_{\;LA}\,.\vspace{-2mm}\]
\end{Thm}\noindent
By Proposition\,3, it follows that  $M(H)$ is not strongly Arens irregular.
\begin{proof}
It is easy to see that $Z_t(\,M(H)'')=
\langle\delta_0\rangle\oplus Z_t(M_0(H)'')$ and that the left annihilator
of the bidual is always part of the (left) topological centre.
Now take $\fr m\in Z^\infty(H)$
(defined in Theorem\,1). We want to show that $\fr m\in Z_t(M_0(H)'')$\,.
\\[1mm]
The property $\fr m\odot L^\infty(H)\subseteq C_0(H)$ means more exactly
$\fr m\odot_a L^\infty(H)\subseteq q(C_0(H))$. This implies (since $q$ is
injective on $C_b(]0,\infty[)$) that $j_d'(\fr m\odot_0f)\in C_0(H)$ for
$f\in L^\infty(H)$ and then $\fr m\odot_0f=j_B(j_d'(\fr m\odot_0f))$ implies
$\fr n\squ\fr m=0$ for $\fr n\in M_{\rm inf}$\,. As noted in the proof of
Corollary\,2 $\fr m\squ L_0^\infty(H)^\perp=\{0\}$ implies $S_1''\fr m=0$\,,
hence by Lemma\,1(a) (extended version) $\fr m\squ\fr n=0$ for
$\fr n\in M_{\rm inf}$\,. The cases where $\fr n\in M_0(H)''$ is carried by a
compact subset of $]0,\infty[$ and where $\fr n$ is carried by $]0,b]$ for
every $b>0$ can be treated as in the first part of the proof of
Theorem\,1.
\\[1mm plus .5mm]
$M_0(H)\subseteq M_{\rm fin}\,,\ Z^\infty(H)
\subseteq M_{\rm inf}$ gives $M_0(H)\cap Z^\infty(H)=\{0\}$\,.
Assume that
$\fr m=\mu+\fr m_1\in(M_0(H)\oplus Z^\infty(H))\cap M_0(H)''_{\;LA}$\,.
Since $M_{\rm fin}$ is a
subalgebra, $M_{\rm inf}$ an ideal, it follows that $\mu\star\nu=0$
for all $\nu\in M_0(H)$, hence $\mu=0$\,.
If $(u_n)$ is an approximate identity  for $L^1(H)$ we get
$\fr m_1=\lim\,\fr m_1\squ u_n=0$\,. Thus we have a direct sum and
$M_0(H)\oplus Z^\infty(H)\oplus M_0(H)''_{\;LA}\subseteq Z_t(M_0(H)'')$\,.
\\
For the converse we start with \,$\fr m\in Z_t(M_0(H)'')$ and we may assume
that \linebreak
$\fr m\vert j_B(C_0(H))=0$\,. $Z_t(\,M_0(H)'')$ is an algebra,
$M_0(H)''\squ M_0(H)''\subseteq L^1(H)''$ (Corollary\,(b) to Proposition\,1)
and clearly $Z_t(M_0(H)'')\cap L^1(H)''\subseteq Z_t(L^1(H)'')$\,.
Thus $\fr m\squ M_0(H)\subseteq Z_t(L^1(H)'')$\,. For
$\mu\in M_0(H)\,,\ f\in C_0(H)$ it is easy to see (using also Lemma\,3(a))
that $\mu\odot_0q(f)\in j_B(C_0(H))$. Thus $\fr m\squ\mu$ vanishes on
$j_B(C_0(H))$ and then it follows from Theorem\,1 that
$\fr m\squ M_0(H)\subseteq Z^\infty(H)$\,. If $Z^\infty(H)$ is trivial, this
entails that $\fr m\in M_0(H)''_{\;LA}$ and the proof is finished. But
in the general case some further arguments are necessary.
\\[0mm plus .5mm]
We have $\fr m=\fr m_1+\fr m_2$ with
$\fr m_1\in M_{\rm fin}\,,\ \fr m_2\in M_{\rm inf}$\,.
Next we want to show that $\fr m_1\in M_0(H)''_{\;LA}$\,. By Lemma\,1(a)
$\fr m\squ\fr n=S''\fr m\squp\fr n$ for $\fr n\in L_0^\infty(H)^\perp$ and
$S''\fr m\in L^1(H)''$\,.
Using Neufang's method as in ($\gamma$) of the proof of Theorem\,1, we obtain
$S''\fr m\in L^1(H)$. Since w*--$\lim_{x\to\infty}\,\nu_x=0$
(Proposition\,1(b))
we have $T\circ q(C_0(H))\subseteq C_0(H)$\,, hence $S''\fr m=0$\,. Using also
that $\supp\nu_x\subseteq[-x,x]$ one can see that $S''$ maps $M_{\rm fin}$
resp. $M_{\rm inf}$ to the corresponding ``finite" resp. ``infinite" part
of $M(\R)''$ and we arrive at $S''\fr m_i=0$ for $i=1,2$. Then the extended
version of Lemma\,1(a) gives $\fr m_1\squ M_{\rm inf}=\{0\}$\,.
\\
$\fr m\squ M_0(H)\subseteq Z^\infty(H)$ implies (by centrality of $\fr m$)
that
$\fr m\squ M_0(H)''\subseteq L_0^\infty(H)^\perp$\,. Since $M_{\rm fin}$ is a
subalgebra, $M_{\rm inf}$ an ideal, it follows that
$\fr m_1\squ M_{\rm fin}=\{0\}$. Combined we get $\fr m_1\in M_0(H)''_{\;LA}$\,.
\\
We can now assume that $\fr m_1=0$\,, i.e. $\fr m\in M_{\rm inf}$\,. Take
$f\in L^\infty(H)$. By Lemma\,3(b) there exists $g\in C_b([0,\infty[)$,
such that $j_d'(\fr m\odot_0f)=j_d'(\fr m\odot_0q(g))$\,. If
$\fr m=\text{w*--}\lim\,\mu_\alpha$ for some net $(\mu_\alpha)\subseteq M_0(H)$
then centrality of $\fr m$ implies that
$\fr m\odot_0q(g)=\lim \mu_\alpha\odot_0q(g)$
{\it weakly}, the same for $f$\,. Then Lemma\,3(a) implies that
$j_d'(\fr m\odot_0f)\in C_b([0,\infty[)$ and
$\fr m\odot_0f=j_B(j_d'(\fr m\odot_0f))$ for all $f\in L^\infty(H)$
\,(this is an analogy to \cite{LL}\,Lemma\,2).
\\
Now let $(u_n)$ be an approximate identity for $L^1(H)$ as in the proof
of Theorem\,1 and $\fr e\in L^1(H)''$ some w*--accumulation point. We have
\,$\lim_{n\to\infty}\langle u_n\squ\fr m,f\rangle=
\lim\,\langle u_n,\fr m\odot_af\rangle=
\lim\,\langle u_n, q(j_d'(\fr m\odot_0f)\rangle=j_d'(\fr m\odot_0f)(0)$ \,for
all $f\in L^\infty(H)$. Thus $\fr e\squ\fr m=\lim u_n\squ\fr m$ holds as a
w*--limit and
by the Grothendieck property, this is a weak limit. It follows that
$\fr e\squ\fr m\in Z^\infty(H)\subseteq Z_t(M_0(H)'')$. Furthermore
(recall that
$\fr e\odot_af=f$) \ $(\fr m\squ\fr e)\odot_0f=\fr m\odot_0(\fr e\odot_af)=
\fr m\odot_0f$ implies
$(\fr e\squ\fr m)\squ\fr n=\fr n\squ(\fr m\squ\fr e)=\fr n\squ\fr m=
\fr m\squ\fr n$ \,for all $\fr n\in M_0(H)''$,
hence $\fr m\,-\,\fr e\squ\fr m\in M_0(H)''_{\;LA}$ finishing the proof.
\vspace{-3mm}\end{proof}
\begin{Rem}	
{\bf \ Compact hypergroups.}\\
For the base space $[0,1]$ (other authors use $[0,\pi]$ or $[0,\frac\pi2]$)
the direct analogue are the {\it Achour--Trim\`eche} hypergroups
(\cite{BH}\,3.5.86). Here a symmetry condition for $A$ is required:
$A(x)=A(1-x)$ for $x\in[0,1]$. It follows that $\{0,1\}$ is a subgroup of
$H$, $\delta_x\star\delta_1=\delta_{1-x}$\,. The formulas for the kernels
in \cite{BS} have been extended to this case in \cite{CMS}\,(3.12). Then
Proposition\,1(a),(c) stays true when restricting to $x,y\in\,]0,1[$\,.
Taking now \,$M_0(H)=
\{\mu\in M([0,1])\!:\, \mu(\{0\})=\mu(\{1\})=0\}\ \ (\,\cong\;M(]0,1[)\,)$,
we get an ideal of codimension $2$ and $M(H)$ is the semidirect product
of $M_0(H)$ with the group algebra of $\{0,1\}$. One can proceed further as
for $\Rp$\,, but now there is no infinite part in $M_0(H)''$,
hence $Z^\infty(H)$ is not present (making the argument
quite a bit easier). One has to consider also elements that are carried by
all neighbourhoods of $1$, but they can be treated like those at $0$.
It follows that $M_0(H)''_{\;LA}$
is always non-trivial and $Z_t(\,M(H)'')=M(H)\oplus M_0(H)''_{\;LA}$, thus
$M(H)$ is not strongly Arens irregular.\\
[.6mm]
For the double coset spaces this covers the cases $G=\SO(n)\,,\;K=\SO(n-1)$
with $n\ge3$ (\cite{BH}\,p.\,246).
\\[.3mm]
Further examples are the {\it dual Jacobi} hypergroups (\cite{BS}\,3.5.81).
For $\alpha>\beta$ the symmetry condition above is not satisfied. It follows
(with some computational effort) from the explicit formulas for the kernels
given in \cite{G}, that in this case $\delta_x\star\delta_1$ is absolutely
continuous for $x>0$ and norm continuity of the product holds also at $1$
(but as far as I see, one can in general no longer use the methods of
\cite{BS}
to get bounds for the densities at $1$). Taking here $M_0(H)=M(]0,1])$ the
structure corresponds to that for $\Rp$. The description of the topological
centres is as in the symmetric case above.\\
For the double coset spaces this covers the cases $G=\U(n)\,,\;
K=\U(n-1)\times\U(1)$ with $n\ge2$ (\cite{BH}\,p.\,246).
\\[.5mm]
It was already observed by Neufang (personal communication)
that for a compact hypergroup $K$,
the algebra $L^1(K)$ is always strongly Arens irregular. This can be done
similarly as in the group case. By \cite{BH}\;Th.\,1.3.28, $K$ admits a
Haar measure and $L^1(K)$ is an ideal in $M(K)$ (\cite{BH}\;Property (1.4.3)).
If $K$ is metrizable, one can proceed now as in \cite{LL} L.\,5 and Cor.\,2
for the trivial subgroup (note that for compact groups and hypergroups only
the first and easier part of the proof of \cite{LL}\;L.\,4 is needed). In
the general case a procedure as on \cite{LL}\;p.\,469 can be applied.
\end{Rem}\vspace{-1mm}
\noindent
We close with some applications for the coset spaces.\vspace{-2mm}
\begin{Pro}    
Let $G$ be a locally compact group, $K$ a compact subgroup.\vspace{.5mm}
Then \hspace*{1.5cm}
$Z_t^{(2)}(L^1(G/K)'')=\{\fr m\in L^1(G/K)'':\lambda_K\squ \fr m\in L^1(G/K)\}$
\qquad and
\\[.5mm] \hspace*{1.5cm}
$Z_t^{(2)}(M(G/K)'')=\{\fr m\in M(G/K)'':\lambda_K\squ \fr m\in M(G/K)\}$\,.
\vspace{-2mm}
\end{Pro}\noindent
$Z_t^{(2)}$ denotes the {\it right} topological centre. $\star$ denotes
now (as at the beginning of the paper) standard convolution on $M(G)$ and
its subalgebras.
For $\squ$ one can take the corresponding Arens product on $M(G)''$. It
restricts to the Arens products on the subalgebras $L^1(G)'',\,M(G/K)'',\dots$.
Since $L^1(G/K)$ is an ideal in $M(G/K)$, it follows that $L^1(G/K)''$ is an
ideal in $M(G/K)''$. $\lambda_K$ denotes normalized Haar
measure of $K$. Note that left and right Arens product coincide if one of the
factors belongs to $M(G)$. The condition $\lambda_K\squ \fr m\in M(G/K)$ is
equivalent to $\lambda_K\squ \fr m\in M(G//K)$ \;(since our embeddings amount to
$M(G//K)=\lambda_K\star M(G/K)$\,).
\\[1mm]
Put $M_{K0}=\{\mu\in M(G/K): \lambda_K\star\mu=0\}\,,\,
L_{K0}=M_{K0}\cap L^1(G/K)$. Then $M(G/K)=M(G//K)\oplus M_{K0}\,,
\,L^1(G/K)=L^1(G//K)\oplus L_{K0}$ and it follows that
$Z_t^{(2)}(M(G/K)'')=M(G//K)\oplus M_{K0}''\,,\,Z_t^{(2)}(L^1(G/K)'')=
L^1(G//K)\oplus L_{K0}''$\,.
If $K$ is not normal in $G$ and $G/K$ is infinite then $L_{K0}$ is infinite
dimensional (see Lemma\;4 below) and it follows
that $L^1(G/K)''$ and $M(G/K)''$ are not strongly Arens irregular.
\begin{proof}
$K\setminus G$ will denote the space of right cosets. Passing to the opposite
group, the claim amounts to determine the left centres
$Z_t(L^1(K\setminus G)'')$ and $Z_t(M(K\setminus G)'')$\,. By analogy with
\cite{LN}\,Lemma\,1 we identify
$L^1(K\setminus G)$ with $\lambda_K\star L^1(G)$, and $M(K\setminus G)$ with
$\lambda_K\!\star M(G)$ \,(the left $K$--invariant elements). $\lambda_K$
is a left identity for $M(K\setminus G)''$, giving
$\fr m\squ\fr n=\fr m\squ\lambda_K\squ\fr n$ for
$\fr m,\fr n\in M(K\setminus G)''$.
Hence $\fr m\in Z_t(L^1(K\setminus G)'')$ iff
$\fr m\squ\lambda_K\in Z_t(L^1(K\setminus G)'')$\,. By a corresponding
identification $M(K\setminus G)\star\lambda_K\linebreak=M(G//K)$ and the same
for $L^1$.
Thus $\fr m\squ\lambda_K\in L^1(G//K)''$. Furthermore,
$\fr m\in Z_t(L^1(K\setminus G)'')\cap L^1(G//K)''$ implies
$\fr m\in Z_t(L^1(G)'')$
\,(this is the opposite case to Proposition\,5 below). The same for $M$\,.
\end{proof}
\begin{Lem}    
If $G$ is a locally compact group, $K$ a compact subgroup such that
$L_{K0}$ is finite dimensional, then either $K$ is normal in $G$ or $G/K$ is
finite.
\end{Lem}\noindent
If $K$ is normal in $G$ then $M_{K0}$ is trivial and if $G/K$ is finite
then $M(G/K)$ is finite dimensional.
\begin{proof}
Assume that $L_{K0}$ is finite dimensional. We have
$v\!-\!\lambda_K\!\star\! v \in L_{K0}$ \, for any $v\!\in\! L^1(G/K)$\,.
If $\mu\in M_{K0}$ and $(v_\iota)$ is a net
in $L^1(G/K)$ converging w* (i.e. for $\sigma(M(G/K),C_0(G/K))$\,)
to $\mu$\,, then $(v_\iota-\lambda_K\star v_\iota)$ converges also to $\mu$\,.
It follows that $L_{K0}$ is w*-dense in $M_{K0}$\,. But finite dimensional
subspaces are closed, hence $M_{K0}=L_{K0}$\,. Recall that $KxK=xK$ is
equivalent to $x\in N_G(K)$ (normalizer of $K$ in $G$\,). We work now on $G$
and assume that $K$ is not normal in $G$\,. If $KxK\neq xK$\, it is easy to
see that
$(\delta_x\star\lambda_K-\lambda_K\star\delta_x\star\lambda_K)(xK)>0$\,.
Now absolute continuity implies $\lambda_G(xK)>0$\,, hence $K$ must be open
in $G$\,. Then by compactness, each double coset of $K$ contains only
finitely many cosets of $K$\,. Finite dimensionality implies that there are
only finitely many double cosets outside $N_G(K)$\,, hence there are only
finitely many cosets of $K$ outside $N_G(K)$\,. It follows that $N_G(K)$
has finite index in $G$\,. For $n,n'\in N_G(K)$ with $nK\neq n'K$ and $x\in G$
one has $xnK\neq xn'K$\,. It follows that $K$ has also finite index in
$N_G(K)$ and combined we get that $G/K$ must be finite.
\end{proof}
\begin{Pro}    
Let $(M,\star)$ be a Banach algebra, $e\in M$ a right identity.
Then\linebreak $Z_t(\,(e\star M)'')=Z_t(M'')\cap (e\star M)''$.\vspace{-2mm}
\end{Pro}
\begin{proof} Easy.\vspace{-2mm}
\end{proof}
\begin{Cor}[\bf to Proposition\,5]\	    
Let $G$ be a locally compact group, $K$ a compact subgroup.
Then \ $Z_t(L^1(G/K)'')\cap L^1(G//K)''=Z_t(L^1(G//K)'')$\qquad and\\
\hspace*{2.9cm}
$Z_t(M(G/K)'')\cap M(G//K)''=Z_t(M(G//K)'')$\,.
\end{Cor}
\begin{Rem}	
\ It follows from Proposition\,4 and the Corollary to
Proposition\,5 that $Z_t(M(G/K)'')\neq Z_t^{(2)}(M(G/K)'')$ whenever
$G//K$ is a Ch\'ebli-\!Trim\`eche hypergroup. Under the condition of
Corollary\,1, one has also
$Z_t(L^1(G/K)'')\neq\linebreak Z_t^{(2)}(L^1(G/K)'')$. Another Example
where the left and right topological centres are different is given in
\cite{LU}\,Example\,2.5 (see also Proposition\,2.10). Further examples and
references are in \cite{DaLa}\,p.\,41, Sections\,6 and 10.
\end{Rem}
\begin{Rem}	
\ In \cite{K} (see in particular the proof of its Lemma\,2.7) it
is shown that if $H$ is any hypergroup admitting a Haar measure (in particular,
if $H$ is commutative or compact or discrete) and $\fr m\in Z_t(L^1(H)'')$ with
$\fr m\ge0$ then $\fr m\in L^1(H)$. In fact, the author claims that this
describes the centre completely, but he does not tell why he believed that
the centre is generated by its non-negative elements.\\
\cite{RKV} asserts some results about the right topological centre of
$L^1(G/K)''$\linebreak ($K$ compact).
Unfortunately, Lemma\,2.2 of \cite{RKV} is wrong,
(1.3) in \cite{RKV} shows a severe misunderstanding about the convolution
in $M(G/K)$\,: \ For $x,y\in G$ the definition (based on the embedding into
$M(G)$, see the Introduction) is
\,$\delta_{xK}\star\delta_{yK}=\int\delta_{xkyK}\,d\lambda_K(k)$.
Unless $y$ normalizes $K$\,, this is always different from $\delta_{xyK}$\,.
\end{Rem}
\begin{Rem}	
\  There is a misleading formulation in \cite{LN} Final
Remark\;(a) on p.\,736. The extension method of Arens can be applied to any
bounded bilinear mapping beween normed spaces (\cite{Da}\;Def.\,A.3.51),
in particular to Banach modules. Then there are corresponding notions of
topological centres and strong Arens irregularity. Using the techniques of
\cite{LN}, one can show that for any locally compact group $G$ with a
compact subgroup $K$\,, one has $Z_t^{(2)}(M(G/K)'')=M(G/K)$ when
considering the second toplogical centre for the left $M(G)$-module
$M(G/K)$. But Proposition\,4 shows that if $K$ is not normal and $G/K$ is
infinite then $Z_t^{(2)}(M(G/K)'')$ gets larger if the toplogical centre
is considered for the Banach algebra $M(G/K)$.
\end{Rem}\vspace{0mm plus 2mm}
\section{Product in Ch\'ebli-\!Trim\`eche hypergroups}\vspace{.5mm}
\label{Proofs}
\begin{proof}[\bf Proof of Proposition 1]\vspace*{0mm plus 5mm}
\item[\bf (a)] Since $\delta_x\star\delta_y=\delta_y\star\delta_x$ it will
be enough to assume $0<y\le x$\,. We consider the differential expression
$\ell(u)(x,y)=u_{xx}(x,y)+\frac{A'(x)}{A(x)}u_x(x,y)
-u_{yy}(x,y)-\frac{A'(y)}{A(y)}u_y(x,y)$\,.\linebreak
By \cite{BH} p.\,209,
\,$\delta_x\star\delta_y$ is related to the solutions of boundary value
problems for $\ell$ \,(the choice of the sign there is not compatible with
other formulas in \cite{BH}, compare \cite{Z}\,(3.1)\,). Put
$\delta_x\star\delta_y=k_{xy}\,\lambda$\,.
For \,$0<a<b<\infty$\,, we will use estimates of Braaksma and de\,Snoo
\cite{BS} to
show uniform integrability of $\{k_{xy}:a\le x,y\le b\}$.
Their setting is slightly different (proving thereby also
absolute continuity of $\delta_x\star\delta_y$), using related
kernels $\tilde k_{xy}(t)$ \,(denoted there as $w_Q(x,y,t)$\,).
Put (for $\alpha_0$ see the Appendix)
\,$\beta(x)=\frac{A'(x)}{A(x)}-\frac{\alpha_0}x\,,\
B(x)=\exp(\frac12\int_0^x\beta(t)\,dt)\,,\ B(x,y)=B(x)\,B(y)$\,.\linebreak
The transmutation $v=B\,u$ gives (see e.g.\;\cite{BH}\,3.5.32)
\,$B\,\ell(u)=\tilde\ell(v)$ with\linebreak
$\tilde\ell(v)(x,y)=
v_{xx}(x,y)+\frac{\alpha_0}x\,v_x(x,y)
-v_{yy}(x,y)-\frac{\alpha_0}y\,v_y(x,y)-q(x,y)\,v(x,y)\;,\
q(x,y)=q(x)-q(y)\,,\
q(x)=\frac{\beta'(x)}2+\frac{\beta^2(x)}4+\frac{\beta(x)\alpha_0}{2x}$\;.
$\tilde\ell$ has the form used in \cite{BS} \,(they write $2p+1$ for our
$\alpha_0$). For a compactly supported $C^2$--function
$\tilde f$ on $]0,\infty[$ a solution of the hyperbolic equation
$\tilde\ell(v)=0$ with boundary condition $v(x,0)=\tilde f(x)$ for
$x>0$ can be obtained by
$v(x,y)=\int_{x-y}^{x+y}\tilde k_{xy}(t)\,\tilde f(t)\,dt$
for $0<y\le x$\,. Going back, write $\tilde f=B\,f$\,, then a solution
for the problem $\ell(u)=0$ with $u(x,0)=f(x)$ is given by
$u(x,y)=\langle\delta_x\star\delta_y\,,f\rangle$ (\cite{BH}\,3.5.19) and
we get $k_{xy}(t)=\frac{B(t)}{B(x)B(y)}\tilde k_{xy}(t)$ for
$0<y\le x\,,\ x-y<t<x+y$\,.
\\[0mm plus1mm]
From $A\ \;C^2$ together with the classical conditions on $A$\,, it follows
that $\beta'$ is bounded (hence $\beta$ is bounded too) on
$]0,c]$ for every $c>0$ (see also the Appendix). Then $x\,q(x)$ is also
bounded on $]0,c]$ and the assumptions in \cite{BS}\,Th.\,2 are satified in
all cases. It relates $\tilde k_{xy}$ to the kernels (denoted as $w_0$) for
the Bessel--Kingman hypergroups (where $A(x)=x^{\alpha_0}\,,\ q=\beta=0$\,).
For $\alpha_0\ge2$ \cite{BS}\,(3.9) and the explicit formula for $w_0$ in
\cite{BS}\,(3.1) (or \cite{BH}\,3.5.61) give that for
$0<a\le x\le b\,,\ 0<y\le x\,,\ x-y<t<x+y$ we have
$\tilde k_{xy}(t)\le \frac cy$
\,for some $c$ depending on $a,b,A$\,. For $0<\alpha_0<2$ the additional term
$\tilde w$ estimated in \cite{BS}\,(3.10)-(3.12) has smaller magnitude
than $w_0$ and one gets
$\tilde k_{xy}(t)\le c\,y^{1-\alpha_0}(y^2-(x-t)^2)^{\frac{\alpha_0}2-1}$ for
$x,y,t$ as above. The factor, when passing to $k_{xy}$\,, is bounded and it
follows that a corresponding estimate holds for $k_{xy}$\,. When adding the
condition $y\ge a$ it follows that the family $\{k_{xy}:a\le x\le y\le b\}$ is
uniformly integrable (recall that $k_{xy}$ is carried by $[x-y,x+y]$\,,
also commutativity of $H$ (\cite{BH}\,Cor.\,3.4.3) implies $k_{xy}=k_{yx}$).
\\[0mm plus1mm]
In \cite{BS} (3.4) an integral formula is given for $\tilde k_{xy}(\xi)$.
The Riemann function and $w_0$ are continuous in $(x,y)$. Using the estimates
made on p.\,69 of \cite{BS} for the integrand (in the various cases) one
can give a common majorant for the integrand when $\xi$ is fixed and $(x,y)$
varies in a sufficiently small neighbourhood
of a given point $(x_0,y_0)$ (the infinities of the integrand are located
at the lines $s+t=\xi$ and $s-t=\xi$\,, these are not affected by $(x,y)$\,).
Then it follows from Lebesgue's convergence
theorem that $\tilde k_{xy}(t)$ and hence also $k_{xy}(t)$ depends
continuously on $(x,y)$\,, i.e. for a sequence $(x_n,y_n)\to(x,y)$
(with $0<y\le x\,,\ y_n\le x_n$) it follows that $k_{x_ny_n}(t)\to k_{xy}(t)$
for each fixed $t$ with $x-y<t<x+y$\,.
Clearly $k_{x_ny_n}(t)\to0$ when
$t<x-y$ or $t>x+y$\,. Now Vitali's convergence theorem implies
$L^1$--convergence $k_{x_ny_n}\to k_{xy}$\,.\vspace{1.2mm plus 2mm}
\item[\bf (c)] Absolute continuity of $\delta_x\star\delta_y$ (with respect to
$\lambda$) implies that $T_yf(x)=\langle\delta_x\star\delta_y\,,f\rangle$ is
defined for $x,y>0$ when $f$ is a bounded (or locally bounded) Lebesgue
measurable function and $T_yf$ depends only on the equivalence class
\!$\mod\lambda$\,. Norm continuity of $x\mapsto\delta_x\star\delta_y$ implies
that $T_yf$ is continuous for $x>0$.
For $0<a\le x\le b$ consider the maximal function
$T^*f(x)=\sup_{0<y\le a}\lvert T_yf(x)\rvert$\,. Put
$\alpha=\min(\frac{\alpha_0}2,1)$\,. It follows from the estimates for
$k_{xy}(t)$ given in (a) that $T^*f(x)\le c\,{\fr M}^\alpha f(x)$ for
the maximal function ${\fr M}^\alpha$ for spherical averages considered
in \cite{SW}\,p.\,1270. By a result of Stein (\cite{SW}\,Th.\,14),
${\fr M}^\alpha$ is $L^p$-continuous for $p>\frac1\alpha$\,. Thus
$T_yf(x)\to f(x)$ for $y\to0+$ pointwise a.e. when $f$ is bounded (or locally
in $L^p$ for some $p>\frac1\alpha$).
\\[1.2mm plus 2mm]
For the Naimark hypergroup (where $A(x)=\sinh^2(x)$ and $\alpha_0=2$)
and more
generally when $\alpha_0\ge2$ the functions $k_{xy}$ are bounded, $T^*$ is
estimated (locally) by the classical Hardy--Littlewood maximal function and
pointwise convergence a.e.\;(even for $f$ locally integrable on $]0,\infty[$\,)
follows from the Lebesgue differentiation theorem.\vspace{.5mm}
\item[\bf (b)] Put $\rho=\frac12\,\lim_{x\to\infty}\frac{A'(x)}{A(x)}$\,,
called the index of the hypergroup $H$ (\cite{BH}\,3.5.47).
For the asymptotics when $x$ is large, we consider
$\ell_{\infty}(u)(x,y)=u_{xx}(x,y)+2\rho u_x(x,y)
-u_{yy}(x,y)-\frac{A'(y)}{A(y)}u_y(x,y)$ \,(\cite{C}\,p.\,451) which can
be defined for $x\in\R\,\ y>0$\,. Solutions for corresponding boundary
value problems can be obtained in a similar way as for $\ell$\,. Since
$\ell_{\infty}$ commutes with translations in the first coordinate, a
corresponding property holds for the measures generating the solutions.
Thus there exist probability measures $\nu_y$ carried by $[-y,y]$ ($y\ge0$),
such that for a $C^2$--function $f$ on $\R$\,, the function
$u(x,y)=(\nu_y\odp f)(x)=\langle\delta_x\starp\nu_y\,,f\rangle=
\int f(x+t)\,d\nu(t)$ satisfies $\ell_{\infty}(u)=0\,,\ u(x,0)=f(x)$\,.
In the special case $A(x)=x^{\alpha_0}$ with $\alpha_0>0$ (coming from the
Bessel--Kingman hypergroups) one obtains $\rho=0$ and
$\ell_{\infty}^{\rm BK}(u)(x,y)=
u_{xx}(x,y)-u_{yy}(x,y)-\frac{\alpha_0}y u_y(x,y)$
(leading to the classical Euler--Poisson--Darboux equation). It is
well known that here $\nu_y^{\rm BK}=g_y^{\rm BK}\lambda$ for $y>0$ with
$g_y^{\rm BK}(t)=\frac1{\mathbf{B}(\frac{\alpha_0}2,\frac12)}\frac1y
(1-(\frac ty)^2)^{\frac{\alpha_0}2-1}$ ($\lvert t\rvert<y\,,\ \mathbf{B}$
denoting the Beta function). One can now make a similar procedure as in (a),
relating the general case to $\ell_{\infty}^{\rm BK}$. Put $v=B_{\infty,0}\,u$
with $B_{\infty,0}(x,y)=e^{\rho x}B(y)$. Then
\,$B_{\infty,0}\:\ell_{\infty}(u)=\tilde\ell_{\infty}(v)$ with
$\tilde\ell_{\infty,0}(v)=\ell_{\infty}^{\rm BK}(v)-q_{\infty,0}\,v$ \,and
$q_{\infty,0}(x,y)=\rho^2-q(y)$ \,($B,q$ as in (a)\,). Repeating the
arguments of \cite{BS} it follows that for $y>0,\ \nu_y=g_y\lambda$ is
absolutely continuous and for $0<y\le b\,,\ \lvert t\rvert<y$ one has
$g_y(t)\le c\,y^{1-\alpha_0}(y^2-(x-t)^2)^{\frac{\alpha_0}2-1}$ \
(for some $c$ depending on $b,A$). Continuing as in (a) it follows that
$y\mapsto\nu_y$ is norm--continuous for $y>0$\,.
\\
For the relationship between $\delta_x\star\delta_y$ and $\nu_y$\,, we
consider another transmutation. Put
$\beta_\infty(x)=\frac{A'(x)}{A(x)}-2\rho\,,\
B_\infty(x,y)=B_\infty(x)=\exp(\frac12\int_1^x\beta_\infty(s)\,ds)$\,, where
$x>0$\,.
Then for $x,y>0$ with $v=B_{\infty}\,u$ one gets
\,$B_{\infty}\:\ell(u)=\tilde\ell_{\infty}(v)$ with
$\tilde\ell_{\infty}(v)=\ell_{\infty}(v)-q_{\infty}\,v$\,, where
$q_\infty(x,y)=q_\infty(x)=
\frac{\beta_\infty'(x)}2+\frac{\beta_\infty^2(x)}4+\beta_\infty(x)\,\rho$ \,.
Now similar estimates as in \cite{BS}\,Th.\,2 can be made. Let
$R_\infty$ be the Riemann function and $\tilde k_{xy}^\infty(s)$ the kernels
for $\tilde\ell_{\infty}$. Then for $0<y\le x\,,\ x-y<\xi<x+y$ the
analogue of \cite{BS}\,(3.4) is
\[g_y(\xi-x)-\tilde k_{xy}^\infty(\xi)=\frac12\,\iint_{\Omega(x,y,\xi,0)}\,
g_t(\xi-s)\,q_\infty(s)\,R_\infty(s,t,x,y)\,d(s,t)\;.
\]
with the rectangle
$\Omega(x,y,\xi,0)=\{(s,t):x-y\le s-t\le\xi\le s+t\le x+y\}$\,.
In estimating the integral (and also the Riemann function in
\cite{BS}\,Th.\,1,
similarly in \cite{CMS}), the authors did  not care so much about the
behaviour of $q$ in the first coordinate by just taking the supremum on $s$\,.
But inspection of the proof shows that for $y\le b$ bounded one can use
instead $\int_{x-y}^{x+y}\lvert q_\infty(s)\rvert\,ds$ in the bound
(this makes it also possible to extend the result to cases where
$q$ or $q_\infty$ are replaced by a locally bounded measure). Using
that $\beta_\infty(x)$ decreases monotonically to $0$ for $x\to\infty$\,, it
is easy to see that this integral tends also to $0$ for $x\to\infty$\,,
uniformly for $0<y\le b$\,. Also the factor
$\frac{B_\infty(\xi)}{B_\infty(x)}$ connecting $k_{xy}(\xi)$ and
$\tilde k_{xy}^\infty(\xi)$ tends to $1$ for $x\to\infty$\,. It follows
that $k_{xy}-\delta_x\starp g_y\to 0$ (norm convergence) for $x\to\infty$\,,
uniformly for $y$ in a compact subset of $]0,\infty[$\,.
\\[0mm plus 2mm]
For the asymptotics of $\nu_x$ when $x\to\infty$ we consider the
eigenfunctions of the (one dimensional) differential operator
$Lf=-f''-\frac{A'}A\,f'$.
For $\lambda\in\C\,,\ \phi_\lambda$ is the solution of the equation
$L\phi_\lambda=(\lambda^2+\rho^2)\,\phi_\lambda$ on $]0,\infty[$ with
$\phi_\lambda(0)=1\,,\ \phi_\lambda'(0)=0$ \,(\cite{BH}\,p.\,223).
Put $\tau_x=e^{-\rho t}\nu_x$\,.
Then by \cite{BH}\,Th.\,3.5.58, $\hat\tau_x(\lambda)=\phi_\lambda(x)$ for
$\lambda\in\C$ \,\linebreak(\,$\hat{ }$ \,denotes the real Fourier and
Fourier--Stieltjes transform,
using the version $\hat f(\lambda)=\int_\R\,e^{-i\lambda t}f(t)\,dt$).
For $A$ unbounded, the main step is to see that for
$\lambda\in\R\setminus\{0\}$ we have
$\phi_\lambda(x)\to0$ for $x\to\infty$\,, i.e. $\phi_\lambda\in C_0(\Rp)$.
This case was left out in the listing of \cite{F}\,Th.\,3.9 (in his
notation $\phi_\lambda=\Phi_{i\,\lambda}$).
\\[0mm plus 2mm]
First, the case $\rho>0$\,. One can use \cite{F}\,L.\,3.7
which gives even exponential decay for the solutions of the differential
equation, thus $\phi_\lambda\in C_0(\Rp)$.
By \cite{BS}\,Prop.\,3.5.49(ii) \,$\phi_0$ is decreasing. Hence (observe
that $\lVert\tau_x\rVert=\phi_0(x)$\,), $\{\tau_x:x>0\}$ is bounded in
$M(\R)$\,.
It follows easily that $\tau_x\to0\;(w^*)$ for $x\to\infty$ \,(it is enough to
consider $\langle\tau_x,f\rangle$ for $f\in C_0\cap L^1(\R)$ with
$\hat f\in L^1$).
Since $\{\nu_x\}$ is bounded too, this implies
$\nu_x\to0\;(w^*)$ for $x\to\infty$\,.
\\[.2mm plus 2mm]
A similar argument (based on Liouville--Green approximation) works for
$\rho=0$ when $A$ is unbounded.
By \cite{E}\,Th.\,2.2.1, there exist independent solutions $y_1,y_2$ of the
differential equation $Lf=\lambda^2\,f$ such that
$y_1(x)\sim A(x)^{-\frac12}\,\exp(I(x))$ and
$y_2(x)\sim A(x)^{-\frac12}\,\exp(-I(x))$ for $x\to\infty$\,, where
$I(x)=\int_1^x\bigl(\frac14(\frac{A'(t)}{A(t)})^2-
\lambda^2\bigr)^\frac12\,dt$\,.
Since $\frac{A'(x)}{A(x)}\to0$ for $x\to\infty$ and $\lambda\neq0$\,,
we have $\frac14(\frac{A'(x)}{A(x)})^2-\lambda^2<0$ for $x\ge x_0$ and it
follows that the real part $\Re I(x)=c_0$ is constant for $x\ge x_0$\,.
This gives
$\lvert y_1(x)\rvert\sim A(x)^{-\frac12}\,e^{c_0}$ and 
$\lvert y_2(x)\rvert\sim A(x)^{-\frac12}\,e^{-c_0}$ for $x\to\infty$.
Since $A(x)\to\infty$ when $\rho>0$\,, this proves $\phi_\lambda\in C_0(\Rp)$.
$\nu_x\to0\;(w^*)$ for $x\to\infty$\, follows from this \;(alternatively,
$\nu_x\to0\;(w^*)$ for $x\to\infty$ follows in this case also from
Proposition\,2(a)\,).
\\[.2mm plus 2mm]
If $A$ is
bounded, $\phi_\lambda\notin C_0(\Rp)$ for $\lambda\in\R\setminus\{0\}$.
Using the recursive equations for $\nu_x$ (see the proof of
Proposition\,2(a)\,) one can show that for the dilated measures
$\langle\nu_x',f\rangle=\int f(\frac tx)\,d\nu_x(t)$ 
(thus $\nu_x'\in M([-1,1])$\,)
one has $\nu_x'\to\frac12(\delta_{-1}+\delta_1)\ (w^*)$ for $x\to\infty$
(if $A$ is bounded). This implies that again $\nu_x\to0\;(w^*)$ for
$x\to\infty$\,.
\end{proof}
\begin{Rem}	
\ Since $\delta_x\star\delta_0=\delta_x$ is atomic, one cannot expect
norm continuity at~$0$.\linebreak
Similarly, it is easy (e.g.\;in the Naimark hypergroup)
to construct examples showing that $T_xf$ need not be continuous at $0$ when
$f$ is discontinuous at $0$\,. In Proposition\,1(c) it is clearly enough
to assume that $f$ is bounded on $[0,b]$ for every $b>0$ (and measurable).
For $x>0\,,\ T_x$ is an integral operator
on $]0,\infty[$\,, hence it exerts some smoothing on functions.
For the Naimark hypergroup it is not hard to see that $T_xf$ is even
Lipschitz continuous on every $[a,b]$ ($0<a<b<\infty$) for $f$
bounded measurable. By \cite{BS}\,(1.3.7), $T_x$
extends to $L^1([0,\infty[,\lambda_A)$, but then $T_xf$ will no longer
be continuous, in general. In the neighbouring case of Levitan hypergroups
(\cite{BH}\,3.5.6; where $\delta_x\star\delta_y$ always has a non-zero atomic
part \cite{BH}\,Th.\,3.5.48) one can show in a similar way that the
absolutely continuous part of $\delta_x\star\delta_y$ depends
norm continuously on $x,y$ for $x,y\in[0,\infty[\times[0,\infty[$\,.
\\[1mm plus .2mm]
In \cite{BX} averages
$S_\varepsilon=
\frac1{\lambda_A([0,\varepsilon])}\int_0^\varepsilon T_x\,d\lambda_A(x)$ were
considered. From the estimates for $k_{xy}$ given in the proof of
Proposition\,1(b), it follows that locally (i.e. on $L^1([a,b])$ for
$1<a<b<\infty$) the
corresponding maximal function can be estimated by the classical
Hardy-Littlewood maximal function. Hence $S_\varepsilon f\to f$ pointwise a.e.
when $\varepsilon\to0+$ holds for $f$ locally integrable on $]0,\infty[$\,.
With some technical effort (using the duality theory of the hypergroup) it
is shown in \cite{BX} that under some regularity conditions on $A$\,, the
maximal function for the $S_\varepsilon$ is of weak type $(1,1)$ on
$L^1(\Rp,\lambda_A)$.
\end{Rem}
\begin{proof}[\bf Proof of Proposition 2]
\item[\bf (a)]
As in \cite{BH}\,3.5.34, the solutions of $\ell_\infty(u)=0$ for $y>0$ with
$u_y(x,0)=0$ satisfy an integral equation. This gives also an
($M(\R)$--valued) integral equation for $\nu_y$\,.
\[ 2\,A(y)\,\nu_y=\int_0^y A'(\eta)
\,\bigl(e^{2\rho(\eta-y)}(\frac{A'(\eta)}{A(\eta)}+2\rho)\,\delta_{\eta-y}\,+
\,e^{2\rho(y-\eta)}(\frac{A'(\eta)}{A(\eta)}-2\rho)\,\delta_{y-\eta}\bigr)
\starp\nu_\eta\,d\eta\,.
\]
If $H$ has sub--exponential growth then $\rho=0$
(\cite{BH}\,Prop.\,3.5.55(i)\,) and the formula simplifies to
\,$2\,A(y)\,\nu_y=\int_0^y A'(\eta)
\,(\delta_{\eta-y}+\delta_{y-\eta})\starp\nu_\eta\,d\eta$\,. Iteration
(inserting the recursion for $\nu_\eta$) gives with
$u_1(\eta)=\frac{A'(\frac\eta2)}{2A(\frac\eta2)}
\ \;(\eta>0),\ \check u_1(\eta)=u_1(-\eta)$\vspace{-2mm}
\begin{multline*}
4\,A(y)\,\nu_y=\int_0^y A'(\eta_1)
\ \bigl(\,\log(\frac{A(y)}{A(\eta_1)})\,(\delta_{\eta_1-y}+\delta_{y-\eta_1})
\\
+(u_1c_{[2\eta_1,2y]})\starp\delta_{-y-\eta_1}
+(\check u_1c_{[-2y,-2\eta_1]})\starp\delta_{y+\eta_1}\,\bigr)
\starp\,\nu_{\eta_1}\,d\eta_1\,.
\end{multline*}
Now let $x>0$ be fixed. Since $u_1$ is decreasing, we have for
$0\!<\!\eta_1\le y-\frac x2\linebreak
\lVert(\delta_0-\delta_x)\starp(u_1c_{[2\eta_1,2y]})\rVert_1=
2\,\log(\frac{A(y)}{A(y-\frac x2)})$\,. And $\rho=0$ implies that
$\frac{A(y)}{A(y-\frac x2)}\to1$ for $y\to\infty$\,. It follows that for
$k=1,2,\dots$ there
exist $y_k$ such that for $y\ge y_k$
\[\lVert\nu_y-\delta_x\starp\nu_y\rVert\le\int_0^y \frac{A'(\eta_1)}{2A(y)}
\,\log(\frac{A(y)}{A(\eta_1)})
\,\lVert\nu_{\eta_1}-\delta_x\starp\nu_{\eta_1}\rVert\,d\eta_1\;+
\;\frac1{2^k}\ .\]
Using that $A$ is unbounded, induction gives
$\lVert\nu_y-\delta_x\starp\nu_y\rVert\le2^{-\frac k2}$ for $y\ge y_k'$\,.
\\
Thus $\lim_{y\to\infty}\lVert\nu_y-\delta_x\starp\nu_y\rVert=0$\,.
From the recursion we get
\;$\lVert\nu_y-\nu_{y+x}\rVert\le\linebreak 2\,(\frac{A(y+x)}{A(y)}-1)+
\frac1{A(y)}\int_0^y A'(\eta)
\,\lVert\nu_\eta-\delta_x\starp\nu_\eta\rVert\,d\eta$\,.
Using again that $A$ is unbounded, it follows that
\;$\lim_{y\to\infty}\lVert\nu_y-\nu_{y+x}\rVert=0$\,.\vspace{.3mm}
\item[\bf (b)]
Multiplying by a constant, we may assume that
$\lim_{y\to\infty}A(y)\,e^{-2\rho y}=1$\,. Monotonicity of $\frac{A'}A$
implies that $A'-2\rho A\ge0$\,, hence $e^{-2\rho y}A(y)$ is increasing,
in particular $A(y)\le e^{2\rho y}$ for $y\ge0$. Put
$u_2(\eta)=2\rho\,e^{2\rho\eta}$ for $\eta\le0$ and $0$ else. Then
$\lVert u_2\rVert_1=1$ hence $\delta_0-\frac{u_2}2$ is invertible in
$M(\R)$\,. We define \,$\nu_\infty=(\delta_0-\frac{u_2}2)^{-1}\starp{ }
\linebreak \frac12\int_0^\infty e^{-2\rho\eta}
(A'(\eta)-2\rho A(\eta))\,\delta_{-\eta}\starp\nu_\eta\,d\eta$
\;(observe that
$e^{-2\rho\eta}(A'(\eta)-2\rho A(\eta))=
\frac{d\ }{d\eta}(e^{-2\rho\eta}A(\eta))$ showing convergence of the integral).
It follows that $\nu_\infty$ satisfies\linebreak
\,$\nu_\infty=
\frac12\int_{-\infty}^0 2\rho\,
e^{2\rho\eta}\delta_{\eta}\!\starp \nu_\infty\,d\eta\,+\,
\frac12\int_0^\infty e^{-2\rho\eta}
(A'(\eta)-2\rho A(\eta))\,\delta_{-\eta}\starp\nu_\eta\,d\eta$\,.
Comparison with the integral equation for $\nu_y$ from the beginning of the
proof gives an estimate\vspace{-6mm}
\begin{multline*} \hspace*{18mm}
\lVert\delta_{-y}\starp\nu_y-\nu_\infty\rVert\ \le\
\frac32(1-A(y)e^{-2\rho y})+
\frac12\,\bigl(\,e^{-4\rho y}\,+\\
\int_0^y(e^{2\rho(\eta-2y)}
\lvert A'(\eta)+2\rho A(\eta)-4\rho e^{2\rho\eta}\rvert+
4\rho\,e^{4\rho(\eta-y)}\lVert\delta_{-\eta}\starp\nu_\eta-\nu_\infty\rVert\,)
\,d\eta\,\bigr)\,.
\end{multline*}
We have
$\lim_{\eta\to\infty}\frac1{A(\eta)}
(A'(\eta)+2\rho A(\eta)-4\rho e^{2\rho\eta})=0$\,. It follows that there
exist $y_k$ such that
$\lVert\delta_{-y}\starp\nu_y-\nu_\infty\rVert\le
2\rho\int_0^ye^{4\rho(\eta-y)}
\lVert\delta_{-\eta}\starp\nu_\eta-\nu_\infty\rVert\,d\eta\,+\,
\frac1{2^{k+1}}$ \ for $y\ge y_k\,,\ k=1,2,\dots$. As in (a) this implies
that $\lim_{y\to\infty}\lVert\delta_{-y}\starp\nu_y-\nu_\infty\rVert=0$\,.
\\[1mm]
By \cite{BH}\,Th.\,3.5.58, $\hat\nu_x(\lambda)=\phi_{\lambda+i\rho}(x)$ for
$\lambda\in\C$\,. Using the asymptotics for $\Phi_\lambda$ in
\cite{BG}\,Th.\,1.17 \,(note that the condition $G\in L^1(1,\infty)$ of
\cite{BG} is equivalent to boundedness of $A(y)\,e^{-2\rho y}$ on
$[1,\infty[$) and the
definition of $c(\lambda)$ in \cite{BG}\,(2.1),\,p.\,248, it follows that
$\hat\nu_\infty(\lambda)=c(-\lambda-i\rho)$ holds for $\lambda\in\R$\,.
$c(-\lambda-i\rho)\neq0$ was shown in \cite{BG}\,pp.\,252,253.
\end{proof}
\begin{Rem}	
\ In \cite{F}\,Th.\,4.4 it is shown that
$\nu_x=\text{w*--}\lim\limits_{y\to\infty}
\,\delta_{-y}\!\starp(\delta_x\star\delta_y)$ \;(convergence in distribution)
holds for general Sturm--Liouville hypergroups.\\
For the Naimark hypergroup it follows easily from the formulas given in the
Introduction that $\nu_x=g_x\lambda$ with $g_x(t)=\frac{e^t}{2\,\sinh(x)}$
for $\lvert t\rvert\le x$ ($0$ else). This gives $\nu_\infty=g_\infty\lambda$
with $g_\infty(t)=e^t$ for $t\le0$ ($0$ else). More generally for the
Bessel--Kingman and the Jacobi hypergroups formulas for $\nu_x$ are
given in \cite{F}\,pp.\,34/35.\\
$c(\lambda)$ is the (generalization of the) Harish-Chandra $c$--function.
For Jacobi hypergroups an explicit formula is given in \cite{BG}\,p.\,212.
Representing the dual of the hypergroup as in \cite{BH}\,Th.\,3.5.50, it
was shown in \cite{BX}\,Th.\,3.26, \cite{BG} Cor.\,0.6 that under some
regularity conditions on $A$ the Plancherel measure is given by
$\frac1{\lvert c(t)\rvert^2}\,\lambda$ \,($t\ge0$).
\end{Rem}
\section*{Appendix: Conditions for Ch\'ebli--\!Trim\`eche hypergroups}
\vspace{1mm} \label{Appen}
The requirements in \cite{BH}\,3.5 for the Sturm--Liouville function $A$
of a Ch\'ebli-\!Trim\`eche hypergroup are as follows.
$A\!:\Rp\to\R$ continuous, increasing,
$A(x)>0$ for $x>0$, \;$A\ \,C^1$
and $\frac{A'}A$ decreasing on $]0,\infty[$\,.
Furthermore they assume that $A$ is unbounded and satisfies a
(rather strong) local condition at $0$\,: \;$A(x)=x^{\alpha_0}a(x)$ with
$\alpha_0>0$, where $a$ is $C^\infty$ in a neighbourhood of $0$, $a(0)>0$,
$a^{(2n+1)}(0)=0$ for $n=0,1,\dots$ (in particular $A(0)=0$). This is designed
to cover the standard examples of the Bessel--Kingman and the Jacobi
hypergroups. We rely on these requirements.
\\
To apply directly the results of \cite{BS} (as described in the proof of
our Proposition\,1) we assumed also that $A$ is $C^2$ on $]0,\infty[$\,.
In general,
monotonicity of $\frac{A'}A$ implies that in any case $A'$ is locally of
bounded variation (see also \cite{F}\,L.\,3.11). As mentioned in the proof of
our Proposition\,1(b) the arguments of \cite{BS} leading to their
Th.\,1\,and\,2 can be extended to the case where $q$ is given by a measure
(satisfying corresponding boundedness conditions), hence this should also
be applicable in the case where $A$ is just $C^1$ (or even locally Lipschitz).
\\[1.2mm plus.5mm]
Similarly to apply the results of \cite{E} in the proof of Proposition\,1(b)
one needs that $A$ is $C^2$ (or at least that $A'$ is locally absolutely
continuous).
The question is also treated in  \cite{F}\,Rem.\,3.13, but the
supplementary argument given there is not correct (resp. incomplete).
In the proof of \cite{E}\,Th.\,1.6.1 an additional transmutation is applied
to the system, leading to \cite{E}\,(1.6.10), before applying Levinson's
theorem. The explicit formulas given by \cite{F}\,p.\,29 show that
in case that $\frac{A'}A$ is not absolutely continuous, the same holds for
the function $Q$ appearing in this transformation. Hence this leads out
of the domain of absolutely continuous functions and it becomes insufficient
to use the classical (pointwise) derivative. Working with the distributional
derivative (given by a measure) this amounts to have a generalization
of Levinson's theorem (\cite{E}\,Th.\,1.3.1) to the case where the
pertubation $R$ is given by bounded measures (instead of integrable
functions). This should be possible (solving the integral equation
corresponding to \cite{E}\,(1.4.13) by iteration).
\\[1.2mm plus 2mm]
\cite{BG} even assume that $A$ is $C^\infty$ throughout $\Rp$\,, to get
refined estimates for the functions $\phi_\lambda$\,. But this is not
needed for their Th.\,1.17.
The solutions $\Phi_{-i\lambda}\,,\Phi_{i\lambda}$ in
\cite{BG}\,Th.\,1.17 (constructed under more restrictive conditions; an
earlier version is in \cite{BX}) correspond to $y_+,y_-$
of \cite{F}\,L.\,3.7. \cite{BG} use an
argument of Olver that is worked out in the Appendix of \cite{BG}.
It needs that $A$ is $C^2$ but could also
be extended to the more general assumptions. The argument of \cite{BS}
that $c(\lambda)\neq0$ (which we used in the proof of our Proposition\,2(b)\,)
is based solely on the asymptotics in their Th.\,1.17.
\\[.5mm plus.5mm]
Concerning the local condition for $A$\,, \cite{BH} keep tacit with references
about the existence and uniqueness theorems for hyperbolic equations
(in particular for the singular case) that are needed in their arguments.
Ch\'ebli and Trim\`eche refer in their papers to earlier results of
Lions and Delsarte who consider extensions to $\Rp$\,, working with
even $C^\infty$--functions and using theorems for the corresponding equations
on the upper half plane. Standard existence theorems for the two dimensional
case (\cite{HW}, \cite{Ka}\,p.\,115) use a different normal form and refer
to a wider class of solutions. Consider the directions
$w_+=\binom11,\;w_-=\binom1{-1}$. In our case, the analogue of the class C*
of \cite{HW} are the functions $u$ continuous for $x,y\ge0,\ C^1$ for $x,y>0$
and such that the directional derivatives $u_{w_+w_-}$ exist and are
continuous for $x,y>0$ and $u_{w_+w_-}=u_{w_-w_+}$ \,(recall that when $u$
is $C^2$ we have $u_{w_+w_-}=u_{xx}-u_{yy}$). In this setting \cite{B}
formulates some existence and uniqueness theorems for the singular case.
Transferring to his notation as in the proof of our Propositon\,1,
his assumptions are satisfied when $a$ is $C^2$, $a(0)>0$ (not necessarily
$a'(0)=0$), giving more general examples of hypergroups with the properties
of Propositions\,1,2.
\\[2mm plus.5mm]
In \cite{SGY}\;Th.\,5 an existence theorem for more general ``Sturm-Liouville
hypergroups" is given, assuming just  Condition SL2 of \cite{BH}\,p.\,202
and that $A,A'$ are locally absolutely continuous.

\end{document}